\documentclass[11pt]{amsart}
\usepackage{amsmath}
\usepackage{graphicx}
\usepackage{amsfonts}
\usepackage{amssymb}
\newtheorem{theorem}{Theorem}
\theoremstyle{plain}
\newtheorem{acknowledgement}{Acknowledgement}

\newtheorem{lemma}{Lemma}

\newtheorem{proposition}{Proposition}
\newtheorem{remark}{Remark}

\numberwithin{equation}{section}

\begin{document}
\title[Vlasov Poisson in bounded domains]{Global existence for the Vlasov-Poisson system in bounded domains}
\author{Hyung Ju Hwang}
\address{Department of Mathematics, Pohang University of Science and Technology, Pohang
790-784, Republic of Korea}
\email{hjhwang@postech.edu}
\author{Juan J. L. Velazquez}
\address{Departamento de Matem\'{a}tica Aplicada. Facultad de Ciencias Matem\'{a}ticas.
Universidad Complutense. Madrid 28040, Spain.}
\email{JJ\_Velazquez@mat.ucm.es}
\date{}
\begin{abstract}In this paper we prove global existence for solutions of the Vlasov-Poisson
system in convex bounded domains with specular boundary conditions and with a
prescribed outward electrical field at the boundary.
\end{abstract}
\maketitle

\section{Introduction}

In this paper we study global solutions for Vlasov-Poisson system in a convex
bounded domain $\Omega$ with specular reflection on the boundary:
\begin{align}
f_{t}+v\cdot\nabla_{x}f+\nabla_{x}\phi\cdot\nabla_{v}f  &
=0\;\;\;\;,\;\;\;\;x\in\Omega\subset\mathbb{R}^{3}\;\;\;,\;v\in\mathbb{R}%
^{3}\;\;,\;\;t>0\label{S1E1}\\
\Delta\phi &  =\int_{\mathbb{R}^{3}}fdv\;\;,\;\;x\in\Omega
\;\;,\;\;t>0\label{S1E2}\\
\frac{\partial\phi}{\partial n_{x}}\left(  t,x\right)   &  =h\left(  x\right)
\;\;,\;\;x\in\partial\Omega\;\;,\;\;t>0\label{S1E3}\\
f\left(  0,x,v\right)   &  =f_{0}\left(  x,v\right)  \;\;\;x\in\Omega
\;\;,\;\;v\in\mathbb{R}^{3}\label{S1E4}\\
f\left(  t,x,v\right)   &  =f\left(  t,x,v^{\ast}\right)  \;\;x\in
\Omega\;\;,\;\;v\in\mathbb{R}^{3}\;\;,\;\;t>0\label{S1E5new}%
\end{align}
where $\Omega$ is a convex bounded domain with $C^{5}$ boundary, $n_{x}$
denotes the outer normal to $\partial\Omega$ and
\begin{equation}
f_{0}\left(  x,v\right)  \geq0\label{S1E4a}%
\end{equation}

Here $f\left(  t,x,v\right)  $ denotes the distribution density of electrons,
$\phi\left(  t,x\right)  $ is the electric potential. The function $h$ in
(\ref{S1E3}) will be assumed to be positive and satisfy the following
compatibility condition:
\begin{equation}
\int_{\Omega}f_{0}\left(  x,v\right)  dxdv=\int_{\partial\Omega}h\left(
x\right)  dS_{x}\label{S1E4b}%
\end{equation}

We will also assume that $f_{0}$ is compactly supported in $\bar{\Omega}%
\times\mathbb{R}^{3}$. In (\ref{S1E5new}) and in the rest of the paper we use
the following notation. Given $\left(  x,v\right)  $ in $\partial\Omega
\times\mathbb{R}^{3}$ we define:
\begin{equation}
v^{\ast}\equiv v-2n_{x}\cdot v\label{S2E0a}%
\end{equation}
where from now on $n_{x}$ is the outer normal vector to $\partial\Omega$ at
the point $x.$

The density of particles at a given point $x$ is given by:
\begin{equation}
\rho\left(  t,x\right)  \equiv\int_{\mathbb{R}^{3}}f\left(  t,x,v\right)
dv\label{S1E1a}%
\end{equation}

In the case of $\Omega=\mathbb{R}^{3},$ the solutions of the system
(\ref{S1E1})-(\ref{S1E4}) are globally defined in time for general initial
data, as it was proved in \cite{Pfaffelmoser} as well as in
\cite{LionsPerthame} using different methods. However, in the case of domains
with boundaries the mathematical theory of well-posedness for the solutions of
the Vlasov-Poisson system is not so complete as in the case of the whole
space. It was proved in \cite{YanGuo} that classical solutions for the problem
(\ref{S1E1})-(\ref{S1E5new}) may not exist in general without the
nonnegativity assumption (\ref{S1E4a}) if $\Omega$ is the half-space
$\mathbb{R}_{+}^{3}$. On the other hand, it was also proved in \cite{YanGuo}
that even with the assumption (\ref{S1E4a}) the derivatives of the solutions
of (\ref{S1E1})-(\ref{S1E5new}) cannot be uniformly bounded near the boundary
of $\Omega$ due to the fact that a Lipschitz estimate for the characteristics
in terms of the initial data is not possible.

One of the main technical difficulties that must be considered in order to
solve (\ref{S1E1})-(\ref{S1E5new}), even for short times, is a careful study
of the evolution of the characteristic curves associated to (\ref{S1E1}) that
remain close during their evolution to the so-called singular set, that is
defined as follows:
\begin{equation}
\Gamma\mathbb{=}\left\{  \left(  x,v\right)  \in\Omega\times\mathbb{R}%
^{3}:x\in\partial\Omega,\ v\in T_{x}\partial\Omega\right\} \label{singset}%
\end{equation}
where $T_{x}\partial\Omega\subset\mathbb{R}^{3}$ is the tangent plane to
$\partial\Omega$ at the point $x.$

Notice that the projection of such characteristic curves\ in the domain
$\Omega$ bounces repeatedly at the boundary $\partial\Omega.$

For the case of a half space $\Omega=\mathbb{R}_{+}^{3}$ the global existence
result was shown in \cite{GuoIndiana} if the initial data $f_{0}$ is assumed
to be constant in a neighbourhood of the singular set. The method there is to
adapt the high-moment technique in \cite{LionsPerthame}. Recently, we have
developed in \cite{HV} a new proof of global existence \ modifying
Pfaffelmoser's idea (cf. \cite{Pfaffelmoser}).

In \cite{Hwang} the case of a general convex bounded domain was considered and
solutions of the linear approximate problem of (\ref{S1E1})-(\ref{S1E5new})
were constructed under the assumption that the initial data $f_{0} $ is
constant near the singular set. The absorbing condition was assumed for the
distribution density $f$ at the boundary to obtain global existence of
solution to the full VP system. Global existence results in $\Omega
=B_{R}\left(  0\right)  $ were also obtained in the paper, under the same
assumption, for a class of radially symmetric data that rule out possible
singular behaviors at the origin.

Actually, several of the technical difficulties that arise in the study of the
characteristics near the singular set had been already addressed in
\cite{GuoIndiana} in the particular case $\Omega=\mathbb{R}_{+}^{3}.$ On the
other hand the results in \cite{Hwang} provide techniques for the problem
(\ref{S1E1})-(\ref{S1E5new}) in more general domains. However, the assumption
of $f_{0}$ being constant near the singular set imposes some restrictions on
the initial data, but its main consequence is to make it possible to ignore
the evolution of the characteristic curves that are close to the singular set.

The main contribution of this paper is to show how to adapt Pfaffelmoser's
ideas and to introduce geometric methods to the problem of general bounded
convex domains with curvatures in order to prove global existence in time for
the solutions of (\ref{S1E1})-(\ref{S1E5new}). It turns out that the effect of
the geometry of the domains modifies in a stronger way than that of the
electric field the dynamics of the characteristics. One of the key ideas
consists in approximating the dynamics of the characteristic curves that are
close to the singular set by means of the dynamics of a hamiltonian system
whose trajectories are constrained to the boundary $\partial\Omega.$ This
approach allows us to include easily in the estimates the effects of the
curvature of the domain. We will then be able to adapt the ideas in
\cite{YanGuo}, \cite{Hwang}, \cite{Pfaffelmoser} to prove global existence.

The plan of the paper is the following. In Section 2 we introduce a coordinate
system that makes it easier to study the trajectories near the singular set
and that we will use in the rest of the paper. In this Section we also
introduce a suitable flatness condition for the initial data $f_{0}$ near the
singular set that will make it possible to obtain solvability for the initial
value problem for the VP system. In Section 3 we describe an iterative
procedure that defines a sequence of functions $\left\{  f_{n}\right\}  $
whose limit as $n\rightarrow\infty$ yields the desired global solution of the
problem. In Section 4 we obtain suitable estimates for the solutions of the
so-called ''linear problem''\ that is the system (\ref{S1E1}), (\ref{S1E4}),
(\ref{S1E5new}) with prescribed $\phi$. In Section 5 we prove that the
sequence $\left\{  f_{n}\right\}  $ converges to a limit function that is
defined as long as a suitable functional $Q\left(  t\right)  $ is bounded.
Section 6 contains some standard energy estimates for the solutions. Section 7
provides a proof of the boundedness of the functional $Q\left(  t\right)  $
adapting the ideas of Pfaffelmoser for this problem, in order to deal with the
geometrical complexity of the domain. This concludes the proof of the Theorem.

\section{Preliminary notation and statement of the main result.}

\subsection{A more convenient coordinate systems near the singular set.}

\bigskip

By assumption $\partial\Omega$ is a $C^{5}$ surface and we will parametrized
it locally using a set of coordinates $\left(  \mu_{1},\mu_{2}\right)  .$ Let
us denote as $x_{\parallel}\left(  \mu_{1},\mu_{2}\right)  $ the point of
$\partial\Omega$ characterized by the values of the parameters $\left(
\mu_{1},\mu_{2}\right)  .$ We will denote as $n\left(  \mu_{1},\mu_{2}\right)
$ the outer normal to $\partial\Omega$ at the point $x_{\parallel}\left(
\mu_{1},\mu_{2}\right)  .$

The Implicit Function Theorem shows that for $\delta>0$ sufficiently small we
can parametrize uniquely the set of points $\partial\Omega+B_{\delta}\left(
0\right)  \subset\mathbb{R}^{3}$ by means of the unique values $\left(
\mu_{1},\mu_{2},x_{\perp}\right)  $ solving the equation:
\begin{equation}
x=x_{\parallel}\left(  \mu_{1},\mu_{2}\right)  -x_{\perp}n\left(  \mu_{1}%
,\mu_{2}\right) \label{S3Ech1}%
\end{equation}

Given $x\in\partial\Omega+B_{\delta}\left(  0\right)  $ we will represent any
vector $v\in\mathbb{R}^{3}$ as:
\begin{equation}
v=v_{\parallel}\left(  \mu_{1},\mu_{2}\right)  -v_{\perp}n\left(  \mu_{1}%
,\mu_{2}\right) \label{S3Ech2}%
\end{equation}
where $v_{\parallel}\left(  \mu_{1},\mu_{2}\right)  \in T_{x_{\parallel
}\left(  \mu_{1},\mu_{2}\right)  }\left(  \partial\Omega\right)  ,\,v_{\perp
}\in\mathbb{R}$ and $\left(  \mu_{1},\mu_{2}\right)  $ are as in
(\ref{S3Ech1}). Moreover, we will represent $v_{\parallel}=v_{\parallel
}\left(  \mu_{1},\mu_{2}\right)  $ using the two coordinates $\left(
w_{1},w_{2}\right)  $ defined by means of:
\begin{equation}
v_{||}=w_{1}u_{1}+w_{2}u_{2}\label{S3Ech3}%
\end{equation}
where $\left\{  u_{1},u_{2}\right\}  $ are the basis of $T_{x_{\parallel
}\left(  \mu_{1},\mu_{2}\right)  }\left(  \partial\Omega\right)  $ given by:%

\begin{equation}
u_{i}=\frac{\partial x_{\parallel}\left(  \mu_{1},\mu_{2}\right)  }%
{\partial\mu_{i}}\;\;,\;\;i=1,2\;.\label{S3Ech4}%
\end{equation}

The system of coordinates $\left(  \mu_{1},\mu_{2},x_{\perp},w_{1}%
,w_{2},v_{\perp}\right)  $ will provide a more convenient representation of
the set of points in the phase space $\Omega\times\mathbb{R}^{3}$ that are
close to the singular set $\Gamma$ defined in (\ref{singset}). The form that
the original equation (\ref{S1E1}) takes in this new set of coordinates is
given in the following Lemma:

\begin{lemma}
\label{Lemma1}The equation (\ref{S1E1}) can be rewritten for $\left(
x,v\right)  \in\left[  \partial\Omega+B_{\delta}\left(  0\right)  \right]
\times\mathbb{R}^{3}$, and using the set of coordinates $\left(  \mu_{1}%
,\mu_{2},x_{\perp},w_{1},w_{2},v_{\perp}\right)  $ in the form:
\begin{equation}
\frac{\partial f}{\partial t}+\sum_{i=1}^{2}\frac{w_{i}}{1+k_{i}x_{\perp}%
}\frac{\partial f}{\partial\mu_{i}}+v_{\perp}\frac{\partial f}{\partial
x_{\perp}}+\sum_{i=1}^{2}\sigma_{i}\frac{\partial f}{\partial w_{i}}%
+F\frac{\partial f}{\partial v_{\perp}}=0\label{S2geom}%
\end{equation}
where:
\begin{equation}
\sigma_{i}\equiv E_{i}-\frac{v_{\perp}w_{i}k_{i}}{1+k_{i}x_{\perp}}%
-\sum_{j,\ell=1}^{2}\frac{\Gamma_{j,\ell}^{i}w_{j}w_{\ell}}{1+k_{j}x_{\perp}%
},\ \ F\equiv E_{\perp}+\sum_{j=1}^{2}\frac{w_{j}^{2}b_{j}}{1+k_{j}x_{\perp}%
}\label{S2defs}%
\end{equation}
where $k_{j}$ are the principal curvatures, $b_{j}$\textbf{\ }are the
coefficients $e$ and $g$ from the second fundamental form according to the
notation in \cite{S} and $\Gamma_{j,\ell}^{i}$ are the Christoffel symbols of
the surface $\partial\Omega.$ The vector $E=\nabla_{x}\phi$ has been written
in the form
\begin{equation}
E=E_{1}u_{1}+E_{2}u_{2}-E_{\perp}n\left(  \mu_{1},\mu_{2}\right)
\label{Edecomp}%
\end{equation}
where $u_{1},u_{2}$ are as in (\ref{S3Ech4}).
\end{lemma}

\begin{proof}
The proof of this result is just a standard lengthy change of variables that
makes use of the classical Gauss-Weingarten equations (cf. \cite{S}, page 124).
\end{proof}

\begin{remark}
The choice of coordinates (\ref{S3Ech1}) implies that the coefficient from the
second fundamental form that is denoted as $f$ in \cite{S} is identically zero.
\end{remark}

\begin{remark}
Notice that since the domain $\Omega$ is convex, and due to (\ref{S1E4a}) we
have $F<0.$
\end{remark}

\subsection{Compatibility conditions for the initial data.}

In order to obtain classical solutions of (\ref{S1E1})-(\ref{S1E5new}) we need
to impose some compatibility conditions on the initial data $f_{0}\left(
x,v\right)  $ at the reflection points of $\partial\Omega\times\mathbb{R}^{3}$
(cf. \cite{YanGuo}, \cite{Hwang}). These conditions are the following:
\begin{align}
f_{0}\left(  x,v\right)   &  =f_{0}\left(  x,v^{\ast}\right) \label{S2E0b}\\
v^{\perp}\left[  \nabla_{x}^{\perp}f_{0}\left(  x,v^{\ast}\right)  +\nabla
_{x}^{\perp}f_{0}\left(  x,v\right)  \right]  +2E^{\perp}\left(  0,x\right)
\nabla_{v}^{\perp}f_{0}\left(  x,v\right)   &  =0\label{S2E0c}%
\end{align}
where $E^{\perp}\left(  0,x\right)  $ is the decomposition of the field
$E\left(  0,x\right)  $ given by (\ref{Edecomp}) and $\nabla_{x}^{\perp
},\,\nabla_{v}^{\perp}$ are the normal components to $\partial\Omega$ of the
gradients $\nabla_{x},\;\nabla_{v}$ respectively.

\subsection{Flatness condition.}

\bigskip

The usual way of dealing with the impossibility of obtaining smooth solutions
for general initial data $f_{0}$ near the singular set consists is assuming
that $f_{0}$ is constant near such a set (cf. \cite{GuoIndiana} as well as
\cite{Hwang}). More precisely we will assume that $f_{0}\in C^{1,\mu} $
satisfies the following flatness condition near the singular set $\Gamma$:%

\begin{equation}
f_{0}\left(  x,v\right)  =\text{constant}\;,\;\operatorname*{dist}\left(
\left(  x,v\right)  ,\Gamma\right)  \leq\delta_{0}\label{S2flatness}%
\end{equation}
for some $\delta>0$ small.

We need to introduce some functional spaces for technical reasons. We define
for $\mu\in\left(  0,1\right)  :$%
\begin{align*}
\left\|  f\right\|  _{C^{1,\mu}\left(  \bar{\Omega}\times\mathbb{R}%
^{3}\right)  }  &  =\sup_{\left(  x,v\right)  ,\left(  x^{\prime},v^{\prime
}\right)  \in\bar{\Omega}\times\mathbb{R}^{3}}\left(  \frac{\left|  \nabla
f\left(  x,v\right)  -\nabla f\left(  x^{\prime},v^{\prime}\right)  \right|
}{\left|  x-x^{\prime}\right|  ^{\mu}+\left|  v-v^{\prime}\right|  ^{\mu}%
}\right)  +\left\|  f\right\|  _{L^{\infty}\left(  \bar{\Omega}\times
\mathbb{R}^{3}\right)  }\;\;,\;\;\nabla=\left(  \nabla_{x},\nabla_{v}\right)
\\
C_{0}^{1,\mu}\left(  \bar{\Omega}\times\mathbb{R}^{3}\right)   &  =\left\{
f\in C^{1,\mu}\left(  \bar{\Omega}\times\mathbb{R}^{3}\right)  :f\text{
compactly supported,\ }\left\|  f\right\|  _{C^{1,\mu}\left(  \bar{\Omega
}\times\mathbb{R}^{3}\right)  }<\infty\right\}
\end{align*}%
\begin{align*}
\left\|  f\right\|  _{C_{\;t;\;x}^{1;1,\mu}\left(  \left[  0,T\right]
\times\bar{\Omega}\right)  }  &  \equiv\sup_{x,x^{\prime}\in\bar{\Omega
},\;t,t^{\prime}\in\left[  0,T\right]  }\frac{\left|  \nabla_{x}f\left(
t,x\right)  -\nabla_{x}f\left(  t^{\prime},x^{\prime}\right)  \right|
}{\left|  x-x^{\prime}\right|  ^{\mu}}\\
&  +\left\|  f\right\|  _{C\left(  \left[  0,T\right]  \times\bar{\Omega
}\right)  }+\left\|  f_{t}\right\|  _{C\left(  \left[  0,T\right]  \times
\bar{\Omega}\right)  }%
\end{align*}%
\begin{align*}
&  \left\|  f\right\|  _{C_{t;\left(  x,v\right)  }^{1;1,\mu}\left(  \left[
0,T\right]  \times\Omega\times\mathbb{R}^{3}\right)  }\\
&  \equiv\sup_{x,x^{\prime}\in\bar{\Omega},\;t,t^{\prime}\in\left[
0,T\right]  }\frac{\left|  \nabla_{x}f\left(  t,x,v\right)  -\nabla
_{x}f\left(  t^{\prime},x^{\prime},v\right)  \right|  +\left|  \nabla
_{v}f\left(  t,x,v\right)  -\nabla_{v}f\left(  t^{\prime},x^{\prime}%
,v^{\prime}\right)  \right|  }{\left|  x-x^{\prime}\right|  ^{\mu}+\left|
v-v^{\prime}\right|  ^{\mu}}+\\
&  +\left\|  f\right\|  _{C\left(  \left[  0,T\right]  \times\bar{\Omega
}\times\mathbb{R}^{3}\right)  }+\left\|  f_{t}\right\|  _{C\left(  \left[
0,T\right]  \times\bar{\Omega}\times\mathbb{R}^{3}\right)  }%
\end{align*}

We define the spaces $C\left(  \left[  0,T\right]  \times\bar{\Omega}\right)
,\;C\left(  \left[  0,T\right]  \times\bar{\Omega}\times\mathbb{R}^{3}\right)
$ as the spaces of continuous functions bounded in the uniform norm.

\begin{remark}
Theorem \ref{globalexistence} that is the main result of this paper can be
proved under a more general condition than (\ref{S2flatness}), namely under
the a vanishing condition similar to the one used in (\cite{YanGuo}) in the
half-line case. More precisely, Theorem \ref{globalexistence} is valid under
the assumption:
\begin{equation}
\left|  f_{0}\left(  x,v\right)  \right|  +\frac{1}{x_{\perp}+v_{\perp}^{2}%
}\left|  \frac{\partial f_{0}}{\partial x}\left(  x,v\right)  \right|
+\frac{1}{x_{\perp}+v_{\perp}^{2}}\frac{1}{v_{\perp}}\left|  \frac{\partial
f_{0}}{\partial v}\left(  x,v\right)  \right|  \leq C\left(  x_{\perp
}+v_{\perp}^{2}\right)  ^{\theta}\;,\;\theta>1\label{vanishing}%
\end{equation}
We will explain at the relevant points the modifications that should be needed
in the proof. However, we have focused mostly in the details of the proof
under the more stringent assumption (\ref{S2flatness}) for simplicity.
\end{remark}

\subsection{The main result: Global existence Theorem.}

\bigskip

The main result of this paper is the following:

\begin{theorem}
\label{globalexistence}Let $f_{0}\in C_{0}^{1,\mu}\left(  \Omega
\times\mathbb{R}^{3}\right)  ,$ $f_{0}\geq0$ for some $0<\mu<1$ and satisfy
(\ref{S2flatness}) and suppose that $h\in C^{2,\mu}\left(  \partial
\Omega\right)  $ satisfies (\ref{S1E4b}) and $h\left(  x_{2},x_{3}\right)
>0$. Then there exists a unique solution $f\in C_{t;\left(  x,v\right)
}^{1;1,\lambda}\left(  \left(  0,\infty\right)  \times\Omega\times
\mathbb{R}^{3}\right)  $, $\phi\in C_{t;x}^{1;3,\lambda}\left(  \left[
0,\infty\right)  \times\Omega\right)  ,$ for some $0<\lambda<\mu,$ of the
Vlasov-Poisson system (\ref{S1E1})-(\ref{S1E4a}) with compact support in $x$
and $v.$
\end{theorem}

\begin{remark}
\label{Rem}Theorem \ref{globalexistence} could be derived using similar
arguments for the case that the function $h$ depends on time and that
$\frac{\partial h}{\partial t}$ is smooth enough. The arguments would require
minor changes but we will just give the details of the argument in the case
that $\frac{\partial h}{\partial t}=0$ for simplicity.
\end{remark}

\section{Iterative procedure}

The usual procedure of proving the existence of solutions for Vlasov-Poisson
models consists in obtaining such a solution as the limit of a sequence of
functions $f^{n}$ that are defined by means of an iterative procedure. More
precisely, we define:
\begin{equation}
f^{0}\left(  t,x,v\right)  =f_{0}\left(  x,v\right)  \;\;,\;\;t\geq
0,\;x\in\Omega,\;v\in\mathbb{R}^{3}\label{itf0}%
\end{equation}%
\begin{align}
f_{t}^{n}+v\cdot\nabla_{x}f^{n}+\nabla_{x}\phi^{n-1}\cdot\nabla_{v}f^{n}  &
=0\;\;\;\;,\;\;\;\;x\in\Omega\subset\mathbb{R}^{3}\;\;\;,\;v\in\mathbb{R}%
^{3}\;\;,\;\;t>0\label{itfn1}\\
\Delta\phi^{n-1}  &  =\rho^{n-1}\left(  x\right)  \equiv\int_{\mathbb{R}^{3}%
}f^{n-1}dv\;\;,\;\;x\in\Omega\;\;,\;\;t>0\label{itfn2}\\
\frac{\partial\phi^{n-1}}{\partial n}  &  =h\;\;,\;\;x\in\partial
\Omega\;\;,\;\;t>0\label{itfn3}\\
f^{n}\left(  0,x,v\right)   &  =f_{0}\left(  x,v\right)  \;\;\;x\in
\Omega\;\;,\;\;v\in\mathbb{R}^{3}\label{itfn4}\\
f^{n}\left(  t,x,v\right)   &  =f^{n}\left(  t,x,v^{\ast}\right)
\;\;x\in\partial\Omega\;\;,\;\;v\in\mathbb{R}^{3}\;\;,\;\;t>0\label{itfn5}%
\end{align}
for $n=1,2,....\;$We assume that $f_{0},\;h$ satisfy (\ref{S1E4a}),
(\ref{S1E4b}) as well as (\ref{S2flatness}).

\bigskip We will also use the notation:
\begin{equation}
E^{n}=\nabla\phi^{n}\label{itfield}%
\end{equation}

Our goal is to show that the sequence $f^{n}$ converges as $n\rightarrow
\infty$ for all $0\leq t<\infty.$ To this end we need to show as a first step
that this sequence is globally defined in time for each $n\geq0.$

\section{Linear problem.}

In order to show that the sequence $\left\{  f^{n}\right\}  $ is well defined
we first study the well-posedness of the problem (\ref{S1E1}), (\ref{S1E4}),
(\ref{S1E5new}) under the assumption that the field $E=\nabla\phi$ is given
and satisfies suitable smoothness conditions. Closely related results have
been obtained in \cite{YanGuo} for the half-line case where geometric
complications are not present.

For further reference we define the evolution of the characteristic curves
associated to (\ref{S1E1}), (\ref{S1E5new}). More precisely, given the field
$E=\nabla_{x}\phi$ we define for each $\left(  x,v\right)  \in\Omega
\times\mathbb{R}^{3}$ the generalized characteristic curve $\left(  X\left(
s;t,x,v\right)  ,V\left(  s;t,x,v\right)  \right)  $ by the following
differential equations:
\begin{align}
\frac{dX}{ds}  &  =V\label{S1E6a}\\
\frac{dV}{ds}  &  =E=\nabla_{x}\phi\label{S1E6b}\\
X\left(  t;t,x,v\right)   &  =x\;\;,\;\;V\left(  t;t,x,v\right)
=v\label{S1E6c}%
\end{align}
as long as $X\in\Omega.$ We extend this definition to arbitrarily long times
assuming that at the times $s=s^{\ast}$ when $X_{n}\left(  s^{\ast
};t,x,v\right)  \in\partial\Omega$ the velocity $V$ bounces elastically at the
boundary, i.e.:
\begin{equation}
V\left(  \left(  s^{\ast}\right)  ^{+};t,x,v\right)  \equiv\lim_{s\rightarrow
s^{\ast},\;s>s^{\ast}}V\left(  s;t,x,v\right)  =\left(  V\left(  \left(
s^{\ast}\right)  ^{-};t,x,v\right)  \right)  ^{\ast}\equiv\left(
\lim_{s\rightarrow s^{\ast},\;s<s^{\ast}}V\left(  s;t,x,v\right)  \right)
^{\ast}\label{S1E6d}%
\end{equation}
where $\left(  \cdot\right)  ^{\ast}$ is as in (\ref{S2E0a}).

\begin{theorem}
\label{Thlinear} Assume that $E\in C_{\;t;\;x}^{0;1,\mu}\left(  \left[
0,T\right]  \times\bar{\Omega}\right)  $ for some $\mu\in\left(  0,1\right)
,$ and $E\cdot n=h>0$ at $\partial\Omega.$ Suppose that $f_{0}\in C_{0}%
^{1,\mu}\left(  \bar{\Omega}\times\mathbb{R}^{3}\right)  ,$ $f_{0}\geq0 $ for
some $\mu>0$. Then there exists a unique $f\in$ $C_{t;\left(  x,v\right)
}^{1;1,\lambda}\left(  \left[  0,T\right]  \times\Omega\times\mathbb{R}%
^{3}\right)  ,$ solution to the linear Vlasov-Poisson system (\ref{S1E1}),
(\ref{S1E4}), (\ref{S1E5new}), for some $0<\lambda<\mu.$ Moreover the function
$f$ satisfies:
\begin{align}
f &  \geq0\;\label{S2Xesp1}\\
\int f\left(  t,x,v\right)  dxdv &  =\int f_{0}\left(  x,v\right)
dxdv\;\;,\;\;t\in\left[  0,T\right] \label{S2Xesp2n}%
\end{align}
\end{theorem}

We will introduce a new coordinate system that will be convenient to study the
dynamics of the characteristic curves for bouncing trajectories. Suppose that
$\left(  \mu_{1},\mu_{2},x_{\perp},w_{1},w_{2},v_{\perp}\right)  $ are as in
(\ref{S3Ech1})-(\ref{S3Ech3}). We then define two new coordinates $\left(
\alpha\left(  t,\mu_{1},\mu_{2},x_{\perp},w_{1},w_{2},v_{\perp}\right)
,\;\beta\left(  t,\mu_{1},\mu_{2},x_{\perp},w_{1},w_{2},v_{\perp}\right)
\right)  $ as follows:
\begin{align}
\alpha\left(  t,\mu_{1},\mu_{2},x_{\perp},w_{1},w_{2},v_{\perp}\right)   &
=\frac{v_{\perp}^{2}}{2}-F\left(  t,\mu_{1},\mu_{2},0,w_{1},w_{2}\right)
x_{\perp},\label{S2C1}\\
\beta\left(  t,\mu_{1},\mu_{2},x_{\perp},w_{1},w_{2},v_{\perp}\right)  -2\pi
H\left(  t,\mu_{1},\mu_{2},x_{\perp},w_{1},w_{2},v_{\perp}\right)   &
=\pi\left(  1-\frac{v_{\perp}}{\sqrt{2\alpha}}\right)  .\label{S2C2}%
\end{align}
where the function $H$ will increase by one at each bounce and $F$ is as in
(\ref{S2defs}). Since $v_{\perp}$ changes from $-\sqrt{2\alpha}$ to
$\sqrt{2\alpha}$ in each bounce, it follows that $\beta$ is continuous along
characteristics. Notice that $\beta$ is just a coordinate that indicates the
specific point in the surface $\left\{  \alpha=\text{constant}\right\}  $
where the trajectory lies. It does not have the specific meaning of an angle,
although we have normalized their variation by $2\pi$ between bounces. Its
functional form has been chosen only for convenience. In all the following we
will write for simplicity $F\left(  t,0\right)  $ instead of $F\left(
t,\mu_{1},\mu_{2},0,w_{1},w_{2}\right)  $ and we will drop the dependence of
$H$ on the variables $\mu_{1},\mu_{2},x_{\perp},w_{1},w_{2},v_{\perp}$ if
there is no risk of confusion. We can then write:
\begin{align*}
x_{\perp}  &  =-\frac{\alpha}{F\left(  t,0\right)  }\left[  1-\left(
1-\frac{\beta-2\pi H\left(  t\right)  }{\pi}\right)  ^{2}\right]  \;,\ \\
v_{\perp}  &  =\sqrt{2\alpha}\left(  1-\frac{\beta-2\pi H\left(  t\right)
}{\pi}\right)  \;.\;
\end{align*}

\bigskip

Making the change of variables $\left(  t,\mu_{1},\mu_{2},x_{\perp}%
,w_{1},w_{2},v_{\perp}\right)  \rightarrow\left(  t,\mu_{1},\mu_{2}%
,\alpha,w_{1},w_{2},\beta\right)  $ we transform the system (\ref{S2geom}) as follows:

\begin{lemma}
\label{Festimate}Under the assumptions of Theorem \ref{globalexistence} we
have that $F\left(  t,\mu_{1},\mu_{2},0,w_{1},w_{2}\right)  \leq
-\varepsilon_{0}<0$ in $\partial\Omega$ in any time interval $0\leq t\leq
t^{\ast}$. Moreover, in the coordinate system $\left(  t,\mu_{1},\mu
_{2},\alpha,w_{1},w_{2},\beta\right)  $ the problem (\ref{S2geom}) can be
reformulated in $\left[  \partial\Omega+B_{\delta}\left(  0\right)  \right]
\times\mathbb{R}^{3}$ as:
\begin{align}
&  \frac{\partial f}{\partial t}+\sum_{i=1}^{2}\frac{w_{i}}{1+k_{i}x_{\perp}%
}\frac{\partial f}{\partial\mu_{i}}+\sum_{i=1}^{2}\sigma_{i}\frac{\partial
f}{\partial w_{i}}\nonumber\\
&  +\left[  v_{\perp}\left(  F\left(  t,x_{\perp}\right)  -F\left(
t,0\right)  \right)  -x_{\perp}\left\{  \sum_{i=1}^{2}\left(  \frac{w_{i}%
}{1+k_{i}x_{\perp}}\frac{\partial F\left(  t,0\right)  }{\partial\mu_{i}%
}+\sigma_{i}\frac{\partial F\left(  t,0\right)  }{\partial w_{i}}\right)
\right\}  \right]  \frac{\partial f}{\partial\alpha}\nonumber\\
&  +\left[  -\frac{\pi v_{\perp}^{2}}{\left(  2\alpha\right)  ^{3/2}}F\left(
t,0\right)  +\frac{2\pi F\left(  t,0\right)  F\left(  t,x_{\perp}\right)
x_{\perp}}{\left(  2\alpha\right)  ^{3/2}}\right.  -\nonumber\\
&  -\left.  \frac{\pi x_{\perp}v_{\perp}}{\left(  2\alpha\right)  ^{3/2}%
}\left\{  \sum_{i=1}^{2}\left(  \frac{w_{i}}{1+k_{i}x_{\perp}}\frac{\partial
F\left(  t,0\right)  }{\partial\mu_{i}}+\sigma_{i}\frac{\partial F\left(
t,0\right)  }{\partial w_{i}}\right)  \right\}  \right]  \frac{\partial
f}{\partial\beta}\nonumber\\
&  =0,\label{S3Em}%
\end{align}
\end{lemma}

\begin{proof}
The inequality $F\left(  t,\mu_{1},\mu_{2},0,w_{1},w_{2}\right)
\leq-\varepsilon_{0}<0$ is a consequence of the convexity of $\Omega$ which
implies that the term $\sum_{j=1}^{2}\frac{w_{j}^{2}b_{j}}{1+k_{j}x_{\perp}} $
in (\ref{S2defs}) is negative near $\partial\Omega$. On the other hand the
continuity of the function $h$ in (\ref{S1E3}) implies that $E_{\perp}%
=-h\leq-\varepsilon_{0}<0$ whence the result follows. The derivation of the
equation (\ref{S3Em}) follows from (\ref{S2C1})-(\ref{S2C2}).
\end{proof}

\begin{remark}
We can rewrite (\ref{S3Em}) as:
\begin{align}
&  \frac{\partial f}{\partial t}+\sum_{i=1}^{2}\frac{w_{i}}{1+k_{i}x_{\perp}%
}\frac{\partial f}{\partial\mu_{i}}+\sum_{i=1}^{2}\sigma_{i}\frac{\partial
f}{\partial w_{i}}\nonumber\\
&  +\left[  v_{\perp}\left(  F\left(  t,x_{\perp}\right)  -F\left(
t,0\right)  \right)  -x_{\perp}\left\{  \sum_{i=1}^{2}\left(  \frac{w_{i}%
}{1+k_{i}x_{\perp}}\frac{\partial F\left(  t,0\right)  }{\partial\mu_{i}%
}+\sigma_{i}\frac{\partial F\left(  t,0\right)  }{\partial w_{i}}\right)
\right\}  \right]  \frac{\partial f}{\partial\alpha}\nonumber\\
&  +\left[  -\frac{\pi F\left(  t,0\right)  }{\sqrt{2\alpha}}+\frac{2\pi
F\left(  t,0\right)  \left[  F\left(  t,x_{\perp}\right)  -F\left(
t,0\right)  \right]  x_{\perp}}{\left(  2\alpha\right)  ^{3/2}}\right.
-\nonumber\\
&  \left.  \frac{\pi x_{\perp}v_{\perp}}{\left(  2\alpha\right)  ^{3/2}%
}\left\{  \sum_{i=1}^{2}\left(  \frac{w_{i}}{1+k_{i}x_{\perp}}\frac{\partial
F\left(  t,0\right)  }{\partial\mu_{i}}+\sigma_{i}\frac{\partial F\left(
t,0\right)  }{\partial w_{i}}\right)  \right\}  \right]  \frac{\partial
f}{\partial\beta}\nonumber\\
&  =0,\label{S3Em1}%
\end{align}
\end{remark}

\begin{remark}
Notice that the dynamics of the tangential part to $\partial\Omega$ of the
characteristics at the singular set is given by the equations:
\[
\frac{d\mu_{i}}{dt}=w_{i}\;\;,\;\;\frac{dw_{i}}{dt}=\sigma_{i}%
\]
\end{remark}

\begin{remark}
The compatibility conditions (\ref{S2E0b}), (\ref{S2E0c}) imply that the
original data $f_{0}$ written in the variables $\left(  t,\mu_{1},\mu
_{2},\alpha,w_{1},w_{2},\beta\right)  $ is a $C^{1,\mathbb{\mu}}$ function.
Moreover, the estimate (\ref{S2flatness}) implies that $f_{0}$ is constant for
$0\leq\alpha\leq C\delta_{0}$ for some $C>0$ fixed.
\end{remark}

\subsection{Velocity Lemma.}

The following result has been proved in \cite{YanGuo}, \cite{Hwang} in a
slightly different manner.

\begin{lemma}
\label{alphbound}Given $\delta>0$ fixed we define $\Gamma_{\delta}%
\equiv\left(  \left[  \partial\Omega+B_{\delta}\left(  0\right)  \right]
\cap\Omega\right)  \times\mathbb{R}^{3}.$ \ Suppose that $E$ satisfies the
regularity assumptions in Theorem \ref{Thlinear}. Then, the characteristic
equations (\ref{S1E6a})-(\ref{S1E6c}) can be solved in the interval
$t\in\left[  0,T\right]  $ for any $\left(  x,v\right)  \in\bar{\Omega}%
\times\mathbb{R}^{3}.$ Moreover, there exist positive constants $C_{1}%
,\;C_{2}$ depending only on $T,\;f_{0},\;\left\|  \nabla E\right\|
_{L^{\infty}\left(  \left[  0,T\right]  ,C^{1/2}\left(  \Omega\right)
\right)  }$ such that, for any $\left(  x,v\right)  \in\Gamma_{\delta\text{ }%
}$ the following estimate holds:
\begin{equation}
C_{1}\left(  v_{\perp}^{2}\left(  0\right)  +x_{\perp}\left(  0\right)
\right)  \leq\left(  v_{\perp}^{2}\left(  t\right)  +x_{\perp}\left(
t\right)  \right)  \leq C_{2}\left(  v_{\perp}^{2}\left(  0\right)  +x_{\perp
}\left(  0\right)  \right)  \;\;,\;\;0\leq t\leq T\label{alphaest}%
\end{equation}

\begin{proof}
Using the estimate $F\leq-\varepsilon_{0}<0$ in Lemma \ref{Festimate} as well
as (\ref{S2C1}), it follows that $\alpha$ is equivalent to $v_{\perp}%
^{2}+x_{\perp}.$ Due to the boundedness of $\Omega$ it is enough to prove this
result for small values of $v_{\perp}^{2}\left(  0\right)  +x_{\perp}\left(
0\right)  ,$ i.e. for points that are close to the singular set. Using
(\ref{S3Em1}) we obtain along the characteristic curves:
\[
\frac{d\alpha}{dt}=v_{\perp}\left(  F\left(  t,x_{\perp}\right)  -F\left(
t,0\right)  \right)  -x_{\perp}\left\{  \sum_{i=1}^{2}\left(  \frac{w_{i}%
}{1+k_{i}x_{\perp}}\frac{\partial F\left(  t,0\right)  }{\partial\mu_{i}%
}+\sigma_{i}\frac{\partial F\left(  t,0\right)  }{\partial w_{i}}\right)
\right\}
\]
Notice, however that keeping this term would not change the essence of the
argument if $\frac{\partial h}{\partial t}$ is smooth as indicated in Remark
\ref{Rem}.

Using our regularity assumptions on $E$ we obtain the estimate:
\[
\left|  \frac{d\alpha}{dt}\right|  \leq C\left\|  \nabla E\right\|
_{C^{1/2}\left(  \Omega\right)  }\left|  v_{\perp}\right|  \left(  x_{\perp
}\right)  ^{1/2}+Cx_{\perp}%
\]
where $C$ depends depends only on $h$ and the geometric properties of
$\partial\Omega.$ Since\ $\frac{v_{\perp}^{2}}{2}+x_{\perp}\leq K\alpha$ for
some positive constant $K$ depending on $h,\;\partial\Omega.$ It then follows
that:
\begin{equation}
\left|  \frac{d\alpha}{dt}\right|  \leq C\alpha\label{eqalpha}%
\end{equation}

Therefore:
\[
C_{1}\alpha\left(  0\right)  \leq\alpha\left(  t\right)  \leq C_{2}%
\alpha\left(  0\right)  \;\;,\;\;0\leq t\leq T
\]
for some $C_{1},\;C_{2}$ depending on $f_{0},\;\left\|  \nabla E\right\|
_{C^{1/2}\left(  \Omega\right)  },$ whence (\ref{alphaest}) follows.
\end{proof}
\end{lemma}

\subsection{Well-posedness of the linear problem.}

\bigskip

\begin{proof}
[\textbf{Proof of Theorem \ref{Thlinear}}]The proof of this Theorem just
follows\ by integrating the equation along the characteristics, combined with
the reflecting boundary condition. The existence and uniqueness of solutions
to (\ref{S1E1})-(\ref{S1E5new}) are not immediate due to the bounces of the
characteristics at the boundary $\partial\Omega,$ but the well-posedness of
this evolution has been obtained in \cite{Hwang}. We will define the sequence
of functions $f_{n}$ by means of the evolution of the characteristics
associated to (\ref{S1E6a})-(\ref{S1E6d}). Then:
\begin{equation}
f\left(  t,x,v\right)  =f_{0}\left(  X\left(  0;t,x,v\right)  ,V\left(
0;t,x,v\right)  \right)  \;\label{S1E6e}%
\end{equation}

In order to show that this procedure defines the function $f$ globally in time
we need to show that for each $\left(  x,v\right)  \in\Omega\times
\mathbb{R}^{3}$ the curves defined by means of (\ref{S1E6a})-(\ref{S1E6d})
intersect the boundary $\partial\Omega\times\mathbb{R}^{3}$ at most a finite
number of times, and in particular they never intersect with the singular set.
This fact is a consequence of the Lemma \ref{alphbound}, since for any
trajectory starting at $\Omega\times\mathbb{R}^{3}$ at time $s=t$ we have
$\alpha>0$ at time $s=t,$ and therefore $\alpha$ remains bounded below and
above during the evolution of the trajectory in the interval $s\in\left[
0,t\right]  .$ Moreover, the characteristics starting in the region $\left\{
\alpha\leq C\delta_{0}\right\}  $ where $f_{0}$ is constant, remain in a set
of the form $\left\{  \alpha\leq C_{1}\left(  T\right)  \delta_{0}\right\}  $
for all $0\leq s\leq T,$ and the characteristics starting in the region
$\left\{  \alpha>C\delta_{0}\right\}  $ where $f_{0}$ is not necessarily
constant remain in a set of the form $\left\{  \alpha>C_{2}\left(  T\right)
\delta_{0}\right\}  $ for $0\leq s\leq T.$ In the first set $f=$constant. In
the second one we have that $\left|  \frac{d\beta}{dt}\right|  $ is bounded by
$\frac{C}{\sqrt{\delta_{0}}}$, where the bound on the number of bounces is
uniformly bounded by $\frac{C}{\sqrt{\delta_{0}}}$. Finally we notice that,
since $E\in C_{x}^{1,\mu},$ classical regularity estimates for the solutions
of ODEs show that the functions $X\left(  s;t,x,v\right)  ,\;V\left(
s;t,x,v\right)  $ are $C^{1,\mu}$ with respect to the variables $\left(
x,v\right)  $ as long as there is no bounces. Moreover, if a trajectory
intersects $\partial\Omega\times\mathbb{R}^{3}$ our regularity assumptions on
$\partial\Omega$ imply also $C^{1,\mu}$ regularity with respect to $\left(
x,v\right)  $ for the values of the coordinates $X$ where the function
$X\left(  s;t,x,v\right)  $ intersects $\partial\Omega,$ as well as for the
time $s=s\left(  t,x,v\right)  $ when such intersection takes place (cf.
\cite{GuoIndiana}, \cite{Hwang}). Therefore, for trajectories bouncing a
finite number of times, the functions $X\left(  s;t,x,v\right)  ,\;V\left(
s;t,x,v\right)  $ can be written as the composition of a finite number of
$C^{1,\mu}$ functions with respect to the variables $\left(  x,v\right)  $.
This proves that the function $f$ defined by means of (\ref{S1E6e}) is
H\"{o}lder with respect to $\left(  x,v\right)  $. The uniqueness of the
solution is due to the fact that the solution is uniquely detemined by the
evolution of the characteristic curves. This concludes the proof of the result.
\end{proof}

\begin{remark}
\label{arcane}The proof of this Theorem is where there would be a relevant
difference if the assumption (\ref{vanishing}) had been used instead of
(\ref{S2flatness}). Indeed, if (\ref{vanishing}) had been assumed, the number
of bounces would not be uniformly bounded for $0\leq t\leq T.$ However, it is
possible to argue as in \cite{YanGuo} to show that Theorem \ref{Thlinear}
holds replacing (\ref{S2flatness}) by (\ref{vanishing}). In our setting the
way of proving this would be to rewrite the characteristic equations
associated to (\ref{S3Em1}) as:
\begin{align}
\frac{d\mu_{i}}{dt} &  =h_{1}\left(  \mu_{i},w_{i},\alpha,\beta\sqrt{\alpha
}\right) \label{eq}\\
\frac{dw_{i}}{dt} &  =h_{2}\left(  t,\mu_{i},w_{i},\alpha,\beta\sqrt{\alpha
}\right) \nonumber\\
\frac{d\alpha}{dt} &  =v_{\perp}\left(  F^{n}\left(  t,x_{\perp}\right)
-F^{n}\left(  t,0\right)  \right)  -x_{\perp}h_{3}\left(  \mu_{i},w_{i}%
,\alpha,\beta\sqrt{\alpha}\right) \nonumber\\
\frac{d\beta}{dt} &  =-\frac{\pi F^{n}\left(  t,0\right)  }{\sqrt{2\alpha}%
}+\frac{2\pi F^{n}\left(  t,0\right)  \left(  F^{n}\left(  t,x_{\perp}\right)
-F^{n}\left(  t,0\right)  \right)  x_{\perp}}{\left(  2\alpha\right)
^{\frac{3}{2}}}+h_{4}\left(  \mu_{i},w_{i},\alpha,\beta\sqrt{\alpha}\right)
\nonumber
\end{align}
where the function $h_{k}$ are smooth in their arguments. The Velocity Lemma
implies that, for a given trajectory the order of magnitude of $\alpha$ does
not change in a time interval $0\leq t\leq T.$ Suppose that we consider
trajectories with $\alpha$ of order $R\leq1.$ Using the change of variables
$\alpha=R\tilde{\alpha},\;\beta=\frac{1}{\sqrt{R}}\tilde{\beta}$ we would
transform (\ref{eq}) in a similar system of equations for the variables
$\tilde{\zeta}=\left(  \mu_{i},w_{i},\tilde{\alpha},\tilde{\beta}\right)  $
with the nonlinearities and their derivatives bounded. Classical regularity
theory for this system would imply that, for $\tilde{\alpha}$ of order one
$\left|  \frac{\partial\tilde{\zeta}}{\partial\tilde{\zeta}_{0}}\right|  $
with $\tilde{\zeta}_{0}=\left(  \mu_{i,0},w_{i,0},\tilde{\alpha}_{0}%
,\tilde{\beta}_{0}\right)  $ would be bounded, as well as the H\"{o}lder norms
evaluated at the points with $\tilde{\alpha}$ of order one. Returning to the
original variables $\alpha,\;\beta$ it would then follow that the worse
derivative would be $\left|  \frac{\partial\beta}{\partial\alpha_{0}}\right|
$ that would be bounded as $\frac{1}{\sqrt{R}}.$ In general, the rescaling for
each factor $\alpha$ would include a term of order $\sqrt{R}.$ Due to the
decay of $f_{0},\;\nabla f_{0}$ near the singular set it would then follows
that $\nabla f$ would be bounded like a power law near the singular set.
Estimates for the H\"{o}lder norms of $\nabla f$ would then be obtained using
the holderianity of $f_{0}$ at distances of order one from the singular set,
and using the decay of $\nabla f$ near the singular set, as well as the fact
that $\frac{\left|  \nabla f\left(  x_{1,}v_{1}\right)  -\nabla f\left(
x_{2,}v_{2}\right)  \right|  }{\left|  \left(  x_{1,}v_{1}^{2}\right)
-\left(  x_{2,}v_{2}^{2}\right)  \right|  ^{\gamma}}\leq C\frac{\left(
\alpha\right)  ^{\beta}}{\alpha^{\gamma}}$ where $\alpha$ is the largest one
of the corresponding value associated to $\left(  x_{1,}v_{1}\right)  $ or
$\left(  x_{2,}v_{2}\right)  .$
\end{remark}

\section{On the convergence of the sequence $\left\{  f_{n}\right\}  .$}

\subsection{Representation formula for the solutions of the Poisson equation
with Neumann boundary conditions.}

The following result is standard. We just include it here for further references:

\begin{proposition}
\label{Neumann}Given a bounded domain $\Omega$ in $\mathbb{R}^{3}$ with a
smooth boundary $\partial\Omega$, there exists a Green's function $G\left(
x,y\right)  $ for the Laplacian operator with Neuman boundary conditions:
\begin{align}
\Delta\phi &  =\rho\left(  x\right)  \;,\;x\in\Omega\;\label{Lap1}\\
\frac{\partial\phi}{\partial n} &  =h\;\;,\;x\in\partial\Omega\label{Lap2}%
\end{align}
with the compatibility condition
\begin{equation}
\int_{\Omega}\rho\left(  x\right)  dx=\int_{\partial\Omega}h\left(  x\right)
dS_{x}\label{Lap3}%
\end{equation}
is given by the following representation formula:
\begin{equation}
\phi\left(  x\right)  =\int_{\Omega}G\left(  x,y\right)  \rho\left(  y\right)
dy-\int_{\partial\Omega}G\left(  x,y\right)  h\left(  y\right)  dS_{y}%
\label{Lap4}%
\end{equation}

Any other solution of (\ref{Lap1})-(\ref{Lap3}) is given by (\ref{Lap4}) up to
an additive constant. The function $G\left(  x,y\right)  $ satisfies the
following estimates:
\begin{align}
\left|  G\left(  x,y\right)  \right|   &  \leq\frac{C}{\left|  x-y\right|
}\;\;,\;\;\left|  \nabla_{x}G\left(  x,y\right)  \right|  \leq\frac{C}{\left|
x-y\right|  ^{2}}\;\label{Lap5}\\
\left|  \nabla_{x}^{2}G\left(  x,y\right)  \right|   &  \leq\frac{C}{\left|
x-y\right|  ^{3}}\;,\text{ }x,y\in\bar{\Omega}\nonumber
\end{align}

where $C$ depends only on the domain $\Omega.$
\end{proposition}

\begin{proof}
This result is well known. See for instance \cite{E}.
\end{proof}

\subsection{The iterative sequence $\left\{  f^{n}\right\}  $ is globally
defined in time.}

\bigskip

We will need the following auxiliary Lemma that states that given $f$ bounded
in the space $C_{t;\left(  x,v\right)  }^{1;1,\lambda}\left(  \left[
0,T\right]  \times\Omega\times\mathbb{R}^{3}\right)  ,$ the corresponding
field $E $ defined by $E=\nabla\phi,$ with $\phi$ satisfying (\ref{S1E2}),
(\ref{S1E3}) satisfies the regularity estimates required in the Theorems
\ref{Thlinear}.

\bigskip Given a function $g:\Omega\rightarrow\mathbb{R,}$ we will denote as
$\left[  \cdot\right]  _{0,\lambda;x}$ the seminorm:
\[
\left[  g\right]  _{0,\lambda;x}\equiv\sup_{x,y\in\Omega}\frac{\left|
g\left(  x\right)  -g\left(  y\right)  \right|  }{\left|  x-y\right|
^{\lambda}}%
\]

We define:%

\begin{equation}
Q\left(  t\right)  \equiv\sup\left\{  \left|  v\right|  ~|~\left(  x,v\right)
\in\text{supp~}f\left(  s\right)  ,\text{ \ }0\leq s\leq t\right\}  .\label{Q}%
\end{equation}

We then have the following result

\begin{lemma}
\label{Lereg}Suppose that $f\in$ $C_{t;\left(  x,v\right)  }^{1;1,\lambda
}\left(  \left[  0,T\right]  \times\Omega\times\mathbb{R}^{3}\right)  .$ Then,
the following estimates hold:
\begin{align}
\left|  \rho\left(  t,x\right)  \right|   &  =\left|  \int f\left(
t,x,v\right)  dv\right|  \leq C\left(  T\right)  \left\|  f\right\|
_{C_{t;\left(  x,v\right)  }^{1;1,\lambda}\left(  \left[  0,T\right]
\times\Omega\times\mathbb{R}^{3}\right)  }\;\;,\;\;x\in\Omega\times\left[
0,T\right] \label{N1}\\
\left|  E\left(  t,x\right)  \right|   &  \leq C\left(  T\right)  \left(
\left\|  f\right\|  _{C_{t;\left(  x,v\right)  }^{1;1,\lambda}\left(  \left[
0,T\right]  \times\Omega\times\mathbb{R}^{3}\right)  }+1\right)
\;\;,\;\;x\in\Omega\times\left[  0,T\right] \label{N1b}\\
\left|  F\left(  t,x\right)  \right|   &  \leq C\left(  T\right)  \left(
\left\|  f\right\|  _{C_{t;\left(  x,v\right)  }^{1;1,\lambda}\left(  \left[
0,T\right]  \times\Omega\times\mathbb{R}^{3}\right)  }+1\right)
\;\;,\;\;x\in\partial\Omega\times\left[  0,T\right] \label{N3a}%
\end{align}%
\begin{equation}
\left|  \nabla\rho\left(  t,x\right)  \right|  \leq\int\left|  \partial
_{x}f\left(  t,x,v\right)  \right|  dv\leq C\left(  T\right)  \left\|
f\right\|  _{C_{t;\left(  x,v\right)  }^{1;1,\lambda}\left(  \left[
0,T\right]  \times\Omega\times\mathbb{R}^{3}\right)  }\;\;,\;\;x\in
\Omega\times\left[  0,T\right] \label{N7n}%
\end{equation}%
\begin{equation}
\left|  \nabla E\left(  t,x\right)  \right|  +\left[  \nabla E\left(
t,\cdot\right)  \right]  _{0,\lambda;x}+\left[  \nabla^{2}E\left(
t,\cdot\right)  \right]  _{0,\lambda;x}\leq C\left(  T\right)  \left\|
f\right\|  _{C_{t;\left(  x,v\right)  }^{1;1,\lambda}\left(  \left[
0,T\right]  \times\Omega\times\mathbb{R}^{3}\right)  }\;\;,\;\;x\in
\Omega\times\left[  0,T\right] \label{N8}%
\end{equation}%
\begin{align}
\left|  \rho_{t}\left(  t,x\right)  \right|   &  \leq C\left(  T\right)
\left\|  f\right\|  _{C_{t;\left(  x,v\right)  }^{1;1,\lambda}\left(  \left[
0,T\right]  \times\Omega\times\mathbb{R}^{3}\right)  }\;\;,\;\;x\in
\Omega\times\left[  0,T\right] \label{N9}\\
\left|  E_{t}\left(  t,x\right)  \right|   &  \leq C\left(  T\right)  \left\|
f\right\|  _{C_{t;\left(  x,v\right)  }^{1;1,\lambda}\left(  \left[
0,T\right]  \times\Omega\times\mathbb{R}^{3}\right)  }\;\;,\;\;x\in
\Omega\times\left[  0,T\right] \label{N10}%
\end{align}
where $C\left(  T\right)  >0$ depends on $Q\left(  T\right)  $ and $T.$
\end{lemma}

\begin{proof}
The estimate (\ref{N1}) is a consequence of the boundedness of the support of
$f$ as well as the boundedness of the $L^{\infty}$ norm of $f.$ The inequality
(\ref{N1b}) follows using a standard regularity theory for the Poisson
equation (cf. \cite{GT}, \cite{M}). The estimate (\ref{N3a}) is a consequence
of the definition of $F$ in (\ref{S2defs}) and our regularity assumptions on
$\partial\Omega$ and $h.$ The inequality (\ref{N7n}) is just a consequence of
the regularity properties of $f$ and the boundedness of its support.Then
(\ref{N8}) is a consequence of classical regularity theory for the Poisson
equation. Similarly we can deduce (\ref{N9}) and (\ref{N10}). Thus the proof
is complete.
\end{proof}

\bigskip

\begin{proposition}
Let $\mu,\lambda\in\left(  0,1\right)  $, satisfying $\mu>\lambda.$ Let
$f_{0}\in C_{0}^{1,\mu}\left(  \bar{\Omega}\times\mathbb{R}^{3}\right)  ,$
$f_{0}\geq0$ satisfy (\ref{S2flatness}). Suppose $h\in C^{1,\mu}\left(
\partial\Omega\right)  ,\;h>0.$ Then, the sequence of functions $f^{n}$ is
globally defined for each $x\in\Omega,\;v\in\mathbb{R}^{3}$ and $0\leq
t<\infty.$ Moreover we have $f^{n}\in C_{t;\left(  x,v\right)  }^{1;1,\lambda
}\left(  \left[  0,T\right]  \times\Omega\times\mathbb{R}^{3}\right)  $ for
any $T>0$ and $\left\|  f^{n}\right\|  _{\infty}=\left\|  f_{0}\right\|
_{\infty},\;\int\rho_{n}\left(  x,t\right)  dx=\int f_{0}\left(  t,x,v\right)  dxdv.$
\end{proposition}

\begin{proof}
We argue by induction. If $n=1$ we use the fact that $\left|  \nabla\phi
_{0}\right|  $ is bounded to obtain that $f^{1}$ is supported in the region
where $\left|  v\right|  \leq C\left(  1+t\right)  .$ Then, $\rho_{1}$ is
bounded by $C\left(  1+t\right)  ^{3}.$ Moreover, due to the regularity of
$h,$ we can apply Theorem \ref{Thlinear} to show that $f^{1}$ is bounded in
$C_{t;\left(  x,v\right)  }^{1;1,\lambda}\left(  \left[  0,T\right]
\times\Omega\times\mathbb{R}^{3}\right)  $ for any $T>0.$ Applying then Lemma
\ref{Lereg} it follows that the norm $\left\|  E^{1}\right\|  _{C_{\;t;\;x}%
^{1;1,\mu}\left(  \left[  0,T\right]  \times\Omega\right)  }$ is bounded. We
can then apply Theorem \ref{Thlinear} to prove that $f^{2}$ is well defined in
$C^{1,\lambda}$ for $0\leq t<\infty.$ Moreover, the support of $f^{2}$ would
be contained in the region $\left|  v\right|  \leq Ch_{2}\left(  t\right)  $
where $h_{2}\left(  t\right)  $ is a continuous increasing function in $t$ and
using again Theorem \ref{Thlinear} it follows that $f^{2}$ is bounded in
$C_{t;\left(  x,v\right)  }^{1;1,\lambda}\left(  \left[  0,T\right]
\times\Omega\times\mathbb{R}^{3}\right)  $ for any $T>0.$ Iterating the
argument we obtain that the sequence $f^{n}$ is defined as indicated. Finally,
using the fact that $f^{n}$ just propagates along the characteristics we
obtain the conservation of $\left\|  f^{n}\right\|  _{\infty}.$ The
conservation of the total mass just follows integrating the equation
(\ref{S1E1}) with respect to the variables $x,v$ whence the proposition follows.
\end{proof}

\bigskip

\subsection{The sequence $\left\{  f^{n}\right\}  $ converges to a solution of
the VP system if the sequence $\left\{  Q^{n}\right\}  $ is bounded.}

\bigskip We define the following measure for the maximal velocities reached
for the distribution $f^{n}:$%

\begin{equation}
Q^{n}\left(  t\right)  \equiv\sup\left\{  \left|  v\right|  ~|~\left(
x,v\right)  \in\text{supp~}f^{n}\left(  s\right)  ,\text{ \ }0\leq s\leq
t\right\}  .\label{Qndef}%
\end{equation}

\begin{proposition}
\label{Propconv}Under the assumptions of Theorem \ref{globalexistence},
suppose that $Q^{n}\left(  t\right)  \leq K$ for $n\geq n_{0},\;0\leq t\leq
T.$ Then, $f^{n}\rightarrow f$ in $C_{t;\left(  x,v\right)  }^{\nu;1,\lambda
}\left(  \left[  0,T\right]  \times\Omega\times\mathbb{R}^{3}\right)  $ as
$n\rightarrow\infty$ with $0<\lambda<\mu,\;0<\nu<1$ and where $f\in
C_{t;\left(  x,v\right)  }^{1;1,\lambda}\left(  \left[  0,T\right]
\times\Omega\times\mathbb{R}^{3}\right)  $ is a solution of (\ref{S1E1}%
)-(\ref{S1E5new}).
\end{proposition}

In the Proof of this Proposition we will use some auxiliary Lemmas. The first
one similar to the one in Theorem 6.2, p 309 of \cite{GuoIndiana}.

\bigskip

\begin{lemma}
Suppose the assumptions on Theorem \ref{globalexistence} are satisfied and
that $Q^{n}\left(  t\right)  \leq K$ for $n\geq n_{0},\;0\leq t\leq T.$ Then:
\begin{equation}
\left|  E^{n}\left(  x,t\right)  \right|  \leq C\left(  T\right)
\label{Ebound}%
\end{equation}%
\begin{equation}
\left[  E^{n}\left(  \cdot,t\right)  \right]  _{C^{\gamma}\left(
\Omega\right)  }\leq C\left(  T\right)  \;\;,\;\text{for any}\;\gamma
\in\left(  0,1\right) \label{Edif}%
\end{equation}
for $n\geq n_{0}+1,\;\;0\leq t\leq T,$ where $C\left(  T\right)  $ depends
only on $K,\;\left\|  f_{0}\right\|  _{L^{\infty}\left(  \Omega\times
\mathbb{R}^{3}\right)  },\;T.$
\end{lemma}

\begin{proof}
Estimate (\ref{Ebound}) follows from classical regularity theory for the
Poisson equation as well as the fact that the density $\rho^{n}$ can be
estimated in $L^{\infty}$ in terms of only $Q^{n}\left(  t\right)  $ and the
initial data.

Indeed, the uniform boundedness of $Q^{n}\left(  t\right)  $ implies that we
can estimate $\rho^{n}$ in the interval $0\leq t\leq T$ uniformly in $n.$
Therefore, regularity theory for the Poisson equation implies that $E^{n}$ is
bounded in $W^{1,p}\left(  \Omega\right)  $ for any $1<p<\infty$ uniformly in
$n.$ Classical embedding results then imply that $E^{n}$ is bounded in
$C^{\gamma}\left(  \Omega\right)  $ for any $\gamma\in\left(  0,1\right)  $
and the result follows.
\end{proof}

We now prove the following basic Lemma that shows that the boundedness of
$Q^{n}$ implies the boundedness of $f$ in the norm $C_{t;\left(  x,v\right)
}^{1;1,\lambda}\left(  \left[  0,T\right]  \times\Omega\times\mathbb{R}%
^{3}\right)  .$

\begin{lemma}
\label{ge}Suppose that $Q^{n}\left(  t\right)  \leq K$ for $n\geq
n_{0},\;0\leq t\leq T.$ Then:
\begin{equation}
\left\|  f^{n}\right\|  _{C_{t;\left(  x,v\right)  }^{1;1,\lambda}\left(
\left[  0,T\right]  \times\Omega\times\mathbb{R}^{3}\right)  }\leq C\left(
T\right)  \;\;,\;\;0\leq t\leq T\label{gradestimate}%
\end{equation}
for $\lambda<\mu,$ $n\geq n_{0}+1,\;$where $C\left(  T\right)  $ depends only
on $K\;$and$\;T.$
\end{lemma}

\begin{proof}
Indeed, the estimates (\ref{Ebound}), (\ref{Edif}) imply that $E^{n}$ is
uniformly bounded on the H\"{o}lder spaces $C_{x}^{\gamma}$ for $0<\gamma<1.$
Notice that choosing $\gamma>\frac{1}{2}$ we can derive the estimates in the
velocity Lemma (cf. \ref{alphbound})\ uniformly in $n.$ In particular this
implies that $\alpha\geq C\delta_{0}$ uniformly in $n$ for $t\in\left[
0,T\right]  .$

We can write the characteristic equations for (\ref{S3Em1}) as:
\begin{align}
\frac{d\mu_{i}}{dt} &  =f_{1}\left(  \mu_{i},w_{i},\alpha,\beta\right)
\label{EqW1}\\
\frac{dw_{i}}{dt} &  =f_{2}\left(  t,\mu_{i},w_{i},\alpha,\beta\right)
\nonumber\\
\frac{d\alpha}{dt} &  =v_{\perp}\left(  F^{n}\left(  t,x_{\perp}\right)
-F^{n}\left(  t,0\right)  \right)  -x_{\perp}f_{3}\left(  \mu_{i},w_{i}%
,\alpha,\beta\right) \nonumber\\
\frac{d\beta}{dt} &  =-\frac{\pi F^{n}\left(  t,0\right)  }{\sqrt{2\alpha}%
}+\frac{2\pi F^{n}\left(  t,0\right)  \left(  F^{n}\left(  t,x_{\perp}\right)
-F^{n}\left(  t,0\right)  \right)  x_{\perp}}{\left(  2\alpha\right)
^{\frac{3}{2}}}+f_{4}\left(  \mu_{i},w_{i},\alpha,\beta\right) \nonumber
\end{align}
where the functions $f_{1},\;f_{3},\;f_{4}$ depend only on the regularity
properties of $\partial\Omega$ and the boundary value $h,$ and therefore are
smooth. The function $f_{2}$ depends also on the field $E.$

We can now take the H\"{o}lder derivative of (\ref{EqW1}) with respect to the
initial data $\zeta_{0}=\left(  \mu_{i,0},w_{i,0},\alpha_{0},\beta_{0}\right)
.$ Let us write also $\zeta=\left(  \mu_{i,},w_{i,},\alpha,\beta\right)  $
and:
\[
\left[  g\right]  _{\lambda;\zeta_{0}}=\sup_{\left|  \left(  x_{0}%
,v_{0}\right)  -\left(  x_{0}^{\prime},v_{0}^{\prime}\right)  \right|  \leq
1}\frac{\left|  g\left(  x_{0},v_{0}\right)  -g\left(  x_{0}^{\prime}%
,v_{0}^{\prime}\right)  \right|  }{\left|  \zeta_{0}-\zeta_{0}^{\prime
}\right|  ^{\lambda}}%
\]

Since, $\alpha$ is uniformly bounded below and using also the uniform
boundedness of $Q^{n}$ it follows that for any $\lambda<\mu$ $\left[
\zeta\right]  _{\lambda;\zeta_{0}}$ satisfies an inequality of the form:
\[
\left|  \frac{d}{dt}\left(  \left[  \zeta\right]  _{\lambda;\zeta_{0}}\right)
\right|  \leq C\left[  \zeta\right]  _{\lambda;\zeta_{0}}%
\]

Therefore $\left[  \zeta\right]  _{\lambda;\zeta_{0}}\leq C\left(  T\right)  $
where $C\left(  T\right)  $ depends only on $K$ and $T.$ Using the formula
$f_{0}\left(  x_{0},v_{0}\right)  =f^{n}\left(  t,x,v\right)  ,$ we obtain a
uniform estimate for $f^{n}\left(  t,\cdot,\cdot\right)  $ in $C_{\left(
x,v\right)  }^{\lambda}\left(  \Omega\times\mathbb{R}^{3}\right)  $ for any
$0<\lambda<\mu.$ It then follows that $\rho^{n}$ is uniformly bounded in
$C_{x}^{\lambda}\left(  \Omega\right)  $ and using classical regularity theory
for the Poisson equation it then follows that $E^{n}\left(  t,\cdot\right)  $
is uniformly bounded in $C_{x}^{1,\lambda}\left(  \Omega\right)  .$ Using
again the lower bound for $\alpha$ we can now apply standard regularity
results for ODEs to (\ref{EqW1}) to obtain that $\left[  \frac{\partial\zeta
}{\partial\zeta_{0}}\right]  _{\lambda;\zeta_{0}}$ is uniformly bounded for
$0\leq t\leq T.$ Therefore $f^{n}$ is uniformly bounded in $C_{\left(
x,v\right)  }^{1,\lambda}\left(  \Omega\times\mathbb{R}^{3}\right)  .$ Using
then the equation (\ref{itfn1}) we obtain uniform estimates for $f_{t}^{n}$ in
$C_{\left(  x,v\right)  }^{\lambda}\left(  \Omega\times\mathbb{R}^{3}\right)
.$ This implies (\ref{gradestimate}) and the Theorem follows.
\end{proof}

\begin{remark}
The result in Lemma \ref{ge} can be proved if (\ref{S2flatness}) is replaced
by (\ref{vanishing}) using the rescaling argument in Remark \ref{arcane}.
\end{remark}

The following simplectic property is a standard consequence of the fact that
the evolution of the characteristics curves is hamiltonian.

\begin{lemma}
Moreover, let us denote as $\left(  X\left(  s;t,x,v\right)  ,V\left(
s;t,x,v\right)  \right)  $ the solution of the characteristic equations
(\ref{S1E6a})-(\ref{S1E6c}). For any $s\in\left[  0,T\right]  $ the
transformation
\[
\left(  x,t\right)  \rightarrow\left(  X\left(  s;t,x,v\right)  ,V\left(
s;t,x,v\right)  \right)
\]
is simplectic. In particular:
\begin{equation}
dX\left(  s;t,x,v\right)  dV\left(  s;t,x,v\right)  =dxdv\label{S2E1c}%
\end{equation}
\end{lemma}

We finally conclude the Proof of Proposition \ref{Propconv}.

\begin{proof}
[\textbf{Proof of Proposition \ref{Propconv}}]We recall that the sequence
$f^{n}$ has been defined by the iteration (\ref{itf0})-(\ref{itfield}).

We claim that the sequence is a Cauchy sequence in $L^{1}\left(  \left[
0,T\right]  \times\Omega\times\mathbb{R}^{3}\right)  $. To prove the claim,
let $f^{n+1}$ and $f^{n}$ be consecutive elements of the sequence $\left\{
f^{n}\right\}  $:
\begin{align}
f_{t}^{n+1}+v\cdot\nabla_{x}f^{n+1}+\nabla_{x}\phi^{n}\cdot\nabla_{v}f^{n+1}
&  =0,\label{n+1}\\
f_{t}^{n}+v\cdot\nabla_{x}f^{n}+\nabla_{x}\phi^{n-1}\cdot\nabla_{v}f^{n}  &
=0.\label{n}%
\end{align}
Substracting (\ref{n}) from (\ref{n+1}) yields
\begin{equation}
\left(  f^{n+1}-f^{n}\right)  _{t}+v\cdot\nabla_{x}\left(  f^{n+1}%
-f^{n}\right)  +\nabla_{x}\phi^{n}\cdot\nabla_{v}\left(  f^{n+1}-f^{n}\right)
=\left(  \nabla_{x}\phi^{n-1}-\nabla_{x}\phi^{n}\right)  \cdot\nabla_{v}%
f^{n}.\label{diff}%
\end{equation}
By integrating (\ref{diff}) along the trajectory $\left(  X\left(  s\right)
,V\left(  s\right)  \right)  $ with $X\left(  t\right)  =x$ and $V\left(
t\right)  =v$ we get:
\begin{align}
\left(  f^{n+1}-f^{n}\right)  \left(  x,v,t\right)   &  =\left(  \left.
f^{n+1}\right|  _{t=0}-\left.  f^{n}\right|  _{t=0}\right)  \left(  X\left(
0\right)  ,V\left(  0\right)  \right) \label{diffChar}\\
&  +\int_{0}^{t}\left(  \nabla_{x}\phi^{n-1}-\nabla_{x}\phi^{n}\right)
\left(  s,X\left(  s\right)  \right)  \cdot\nabla_{v}f^{n}\left(  s,X\left(
s\right)  ,V\left(  s\right)  \right)  ds\nonumber\\
&  =\int_{0}^{t}\left(  \nabla_{x}\phi^{n-1}-\nabla_{x}\phi^{n}\right)
\left(  s,X\left(  s\right)  \right)  \cdot\nabla_{v}f^{n}\left(  s,X\left(
s\right)  ,V\left(  s\right)  \right)  ds\nonumber
\end{align}
where:
\[
\frac{dX}{ds}=V\;\;,\;\;\frac{dV}{ds}=\nabla_{x}\phi^{n}\left(  s,X\left(
s\right)  \right)
\]

Applying the representation formula (\ref{Lap4}) and the estimates
(\ref{Lap5}) to compute the difference $\left(  \nabla_{x}\phi^{n-1}%
-\nabla_{x}\phi^{n}\right)  \left(  s,x\right)  $ we obtain:
\[
\left|  \left(  \nabla_{x}\phi^{n-1}-\nabla_{x}\phi^{n}\right)  \left(
s,x\right)  \right|  \leq C\int_{\Omega}\frac{\left|  \rho^{n}\left(
y\right)  -\rho^{n-1}\left(  y\right)  \right|  }{\left|  x-y\right|  ^{2}}dy
\]

Integrating (\ref{diffChar}) over the phase space $\left(  x,v\right)  $ we
obtain, applying Fubini's theorem:
\begin{align*}
&  \left\|  f^{n+1}\left(  t\right)  -f^{n}\left(  t\right)  \right\|
_{L^{1}\left(  \Omega\times\mathbb{R}^{3}\right)  }\\
&  \leq\int_{0}^{t}\iint_{\Omega\times\mathbb{R}^{3}}\left|  \left(
\nabla_{x}\phi^{n-1}-\nabla_{x}\phi^{n}\right)  \left(  s,X\left(  s\right)
\right)  \right|  \left|  \nabla_{v}f^{n}\left(  s,X\left(  s\right)
,V\left(  s\right)  \right)  \right|  dvdxds\\
&  \leq C\int_{0}^{t}\int_{\Omega}K^{n}\left(  y,s\right)  \left|  \rho
^{n}\left(  s,y\right)  -\rho^{n-1}\left(  s,y\right)  \right|  dyds,
\end{align*}
where:
\begin{align*}
K^{n}\left(  y,s\right)   &  \equiv\iint_{\Omega\times\mathbb{R}^{3}}\frac
{1}{\left|  X\left(  s\right)  -y\right|  ^{2}}\left|  \nabla_{v}f^{n}\left(
s,X\left(  s\right)  ,V\left(  s\right)  \right)  \right|  dX\left(  s\right)
dV\left(  s\right) \\
&  =\iint_{\Omega\times\mathbb{R}^{3}}\frac{1}{\left|  x-y\right|  ^{2}%
}\left|  \nabla_{v}f^{n}\left(  s,x,v\right)  \right|  dxdv.
\end{align*}

In the last identity we used the Liouville principle (\ref{S2E1c}). We now
estimate $K^{n}\left(  x,t\right)  $ as follows:
\begin{align*}
K^{n}\left(  x,t\right)   &  =\int_{\left|  y-x\right|  \leq r}\frac{\left\|
\nabla_{v}f^{n}\left(  t,y,\cdot\right)  \right\|  _{L_{v}^{1}}}{\left|
x-y\right|  ^{2}}dy+\int_{\left|  y-x\right|  \geq r}\frac{\left\|  \nabla
_{v}f^{n}\left(  t,y,\cdot\right)  \right\|  _{L_{v}^{1}}}{\left|  x-y\right|
^{2}}dy\\
&  \leq Cr\left\|  \nabla_{v}f^{n}\left(  t,y\right)  \right\|  _{L_{x}%
^{\infty}\left(  L_{v}^{1}\right)  }+\frac{\left\|  \nabla_{v}f^{n}\left(
t,\cdot\right)  \right\|  _{L_{x}^{1}\left(  L_{v}^{1}\right)  }}{r^{2}}.
\end{align*}
We choose $r$ to optimize the right hand side of the above inequality,
namely:
\[
r\left\|  \nabla_{v}f^{n}\left(  t\right)  \right\|  _{L_{x}^{\infty}\left(
L_{v}^{1}\right)  }=\frac{\left\|  \nabla_{v}f^{n}\left(  t\right)  \right\|
_{L_{x}^{1}\left(  L_{v}^{1}\right)  }}{r^{2}}.
\]
where from now on we use by shortedness the notation $L_{x}^{\infty}\left(
L_{v}^{1}\right)  =L^{\infty}\left(  \Omega,L^{1}\left(  \mathbb{R}%
^{3}\right)  \right)  ,\;L_{x}^{1}\left(  L_{v}^{1}\right)  =L^{1}\left(
\Omega\times\mathbb{R}^{3}\right)  .$ Thus we have
\[
K^{n}\left(  t,x\right)  \leq C\left\|  \nabla_{v}f^{n}\left(  t\right)
\right\|  _{L_{x}^{\infty}\left(  L_{v}^{1}\right)  }^{2/3}\left\|  \nabla
_{v}f^{n}\left(  t\right)  \right\|  _{L_{x}^{1}\left(  L_{v}^{1}\right)
}^{1/3}.
\]
Since $f^{n}$ is bounded in $W^{1,\infty}$ and the supports are bounded
uniformly in $n$ (due to the global bound on $Q^{n}\left(  t\right)  ,$ for
$0\leq t\leq T)$ it is easy to see that $\left\|  \nabla_{v}f^{n}\left(
t\right)  \right\|  _{L_{x}^{1}\left(  L_{v}^{1}\right)  }$ and $\left\|
\nabla_{v}f^{n}\left(  t\right)  \right\|  _{L_{x}^{\infty}\left(  L_{v}%
^{1}\right)  }$ are uniformly bounded in $n$ and hence $K^{n}\left(
t,x\right)  $ is uniformly bounded.

Thus we obtain
\begin{equation}
\left\|  f^{n+1}\left(  t\right)  -f^{n}\left(  t\right)  \right\|
_{L^{1}\left(  \Omega\times\mathbb{R}^{3}\right)  }\leq C\left(  T\right)
\int_{0}^{t}\left\|  f^{n}\left(  s\right)  -f^{n-1}\left(  s\right)
\right\|  _{L^{1}\left(  \Omega\times\mathbb{R}^{3}\right)  }ds\label{gronw1}%
\end{equation}
where $C=C\left(  T\right)  $ depends only on $T,Q\left(  T\right)  ,$ and the
initial data. Notice that (\ref{gronw1}) implies by iteration that:
\begin{equation}
\left\|  f^{n+1}\left(  t\right)  -f^{n}\left(  t\right)  \right\|
_{L^{1}\left(  \Omega\times\mathbb{R}^{3}\right)  }\leq C_{1}\theta
^{n}\label{Cauchy1}%
\end{equation}
for some $\theta<1,$ and $0\leq t\leq\varepsilon_{0},$ if $\varepsilon_{0}$ is
sufficiently small depending only on $T$. Then (\ref{gronw1}) implies:
\[
\left\|  f^{n+1}\left(  t\right)  -f^{n}\left(  t\right)  \right\|
_{L^{1}\left(  \Omega\times\mathbb{R}^{3}\right)  }\leq C\left(  T\right)
C_{1}\theta^{n}+C\left(  T\right)  \int_{\varepsilon_{0}}^{t}\left\|
f^{n}\left(  s\right)  -f^{n-1}\left(  s\right)  \right\|  _{L^{1}\left(
\Omega\times\mathbb{R}^{3}\right)  }ds
\]

Therefore we obtain that (\ref{Cauchy1}) is valid for $t\in\left[
\varepsilon_{0},2\varepsilon_{0}\right]  $ changing $C_{1}$ if needed.
Iterating the argument, it then follows that $\left\{  f^{n}\right\}  $ is a
Cauchy sequence in the space $L^{\infty}\left(  \left[  0,T\right]
,L^{1}\left(  \Omega\times\mathbb{R}^{3}\right)  \right)  .$

Once we know that $f^{n}$ is a Cauchy sequence in $L^{1}\left(  \left[
0,T\right]  \times\Omega\times\mathbb{R}^{3}\right)  $, we can show that the
sequence is Cauchy in $C_{t;\left(  x,v\right)  }^{\nu;1,\lambda}\left(
\left[  0,T\right]  \times\Omega\times\mathbb{R}^{3}\right)  $ for any
$0<\lambda<\mu,\;0<\nu<1$ arguing by interpolation. Indeed, using
(\ref{gradestimate}) with $\lambda$ replaced by $\tilde{\lambda}$ satisfying
$\lambda<\tilde{\lambda}<\mu$ and interpolating between $L^{1}\left(  \left[
0,T\right]  \times\Omega\times\mathbb{R}^{3}\right)  $ and $W^{1,\infty
}\left(  \left[  0,T\right]  \times\Omega\times\mathbb{R}^{3}\right)  $ we
obtain that $f^{n}$ is Cauchy in $W^{1,p}\left(  \left[  0,T\right]
\times\Omega\times\mathbb{R}^{3}\right)  $ for any $p>1$. Using Sobolev's
embeddings we can obtain that $f^{n}$ is a Cauchy sequence in $C^{\tilde{\nu}%
}\left(  \left[  0,T\right]  \times\Omega\times\mathbb{R}^{3}\right)  $ for
any $0<\tilde{\nu}<\nu$. Interpolating then in Schauder spaces between
$C^{1,\lambda}$ and $C^{\tilde{\nu}}$ we obtain the desired convergence. In
order to check that $f\in C_{t;\left(  x,v\right)  }^{1;1,\lambda}\left(
\left[  0,T\right]  \times\Omega\times\mathbb{R}^{3}\right)  ,$ we use
(\ref{n}) to obtain:
\[
f^{n}\left(  t\right)  =f_{0}-\int_{0}^{t}\left[  v\cdot\nabla_{x}f^{n}\left(
s\right)  +\nabla_{x}\phi^{n-1}\left(  s\right)  \cdot\nabla_{v}f^{n}\left(
s\right)  \right]  ds
\]
Passing to the limit in this equation as $n\rightarrow\infty$ it follows
that:
\[
f\left(  t\right)  =f_{0}-\int_{0}^{t}\left[  v\cdot\nabla_{x}f\left(
s\right)  +\nabla_{x}\phi\left(  s\right)  \cdot\nabla_{v}f\left(  s\right)
\right]  ds
\]
and this implies the desired differentiability for $f,$ and the Proposition follows.
\end{proof}

\subsection{Convergence of the sequence $Q^{n}$ to $Q$.}

\begin{proposition}
\label{Qconv}Let $Q^{n},\;Q$ be as in (\ref{Qndef}), (\ref{Q}) respectively.
Suppose that $\max\left\{  \sup_{n\geq n_{0}}Q^{n}\left(  t\right)  ,Q\left(
t\right)  \right\}  \leq K$ for $0\leq t\leq T.$ Let us assume also that
$f^{n}\rightarrow f$ in $C_{t;\left(  x,v\right)  }^{\nu;1,\lambda}\left(
\left[  0,T\right]  \times\Omega\times\mathbb{R}^{3}\right)  $ for any
$0<\lambda<\mu,\;0<\nu<1.$ Then $\lim_{n\rightarrow\infty}Q^{n}\left(
t\right)  =Q\left(  t\right)  \;$uniformly on $\left[  0,T\right]  .$
\end{proposition}

\begin{proof}
The characteristics starting in $\alpha\left(  0\right)  \geq C\delta_{0}$
remain during their evolution in the set $\left\{  \alpha\left(  t\right)
\geq C\delta_{0}\right\}  $ due to Lemma \ref{alphbound}. Therefore, these
characteristics remain separated from the singular set and we can estimate
their difference as $n\rightarrow\infty$ as it was made in \cite{HV}. Indeed,
for these characteristics the number of bounces is uniformly bounded in $n$ in
the interval $0\leq t\leq T$.\textbf{\ }Moreover, the times where these
bounces take place for the functions $f_{n}$ converge to the corresponding
times for the bouncing times for the characteristics associated to $f$ and
using the fact that $E^{n}\rightarrow E$ as $n\rightarrow\infty$ we obtain
also that the corresponding characteristic curves associated to $f^{n}$
converge to the ones associated to $f$ between bounces, and due to the
boundedness of the number of such bounces, it follows that the functions
$V^{n}\left(  s,t;x,v\right)  $ converge uniformly to $V\left(
s,t;x,v\right)  $ as $n\rightarrow\infty.$ Since $Q^{n}\left(  t\right)  $ is
the maximum value of $\left|  V^{n}\right|  $ associated to characteristics
$\left(  X^{n}\left(  s,t;x,v\right)  ,V^{n}\left(  s,t;x,v\right)  \right)  $
which at $s=0$ lie in the support of $f_{0}$, it follows that $Q^{n}\left(
t\right)  \rightarrow Q\left(  t\right)  ,$ and the result follows.
\end{proof}

\bigskip

\subsection{Prolongability of uniform estimates for the functions $f^{n}$}

\begin{proposition}
\label{Prolog}Suppose that for some $T\geq0$ there exist $K>0$ and $n_{0}%
\geq0$ such that for any $n\geq n_{0}$ and $0\leq t\leq T$ we have
$Q^{n}\left(  t\right)  \leq K.$ Then, there exists $\varepsilon
_{0}=\varepsilon_{0}\left(  K,\left\|  f_{0}\right\|  _{\infty}\right)  >0$
such that for $0\leq t\leq T+\varepsilon_{0}$ and $n\geq n_{0}$ the following
estimate holds:
\[
Q^{n}\left(  t\right)  \leq2K
\]
\end{proposition}

\bigskip

\begin{proof}
Notice that:
\begin{equation}
\left|  \rho^{n}\right|  =\left|  \int f^{n}dv\right|  \leq\left\|
f_{0}\right\|  _{\infty}\left(  Q^{n}\left(  t\right)  \right)  ^{3}%
\label{rhonest}%
\end{equation}

Using the representation formula (\ref{Lap4}) for the solutions of the Poisson
equation in $\Omega$ with Neumann boundary conditions in Proposition
\ref{Neumann} we obtain:
\[
\;\;\left|  \nabla\phi^{n}\right|  \leq C\left[  \left\|  f_{0}\right\|
_{\infty}\left(  Q^{n}\left(  t\right)  \right)  ^{3}+\int\rho_{n}dx+\left\|
h\right\|  _{C^{1,\mu}}\right]  \;.
\]

On the other hand, using (\ref{S1E6b}) and (\ref{Qndef}) we obtain the
following inequality for $t\geq T$:
\[
Q^{n+1}\left(  t\right)  \leq Q^{n}\left(  T\right)  +C\left\|  f_{0}\right\|
_{\infty}\int_{T}^{t}\left(  Q^{n}\left(  s\right)  \right)  ^{3}ds+C\left(
t-T\right)
\]
where $C>0$ is just a numerical constant independent of $n$ and $Q^{n}.$ Using
the assumption on $Q^{n}$ we obtain:
\[
Q^{n+1}\left(  t\right)  \leq K+C\left\|  f_{0}\right\|  _{\infty}\int_{T}%
^{t}\left(  Q^{n}\left(  s\right)  \right)  ^{3}ds+C\left(  t-T\right)
\;\;,\;\;t\geq T
\]
Defining $R^{n}\left(  t\right)  =\max\left\{  Q^{\ell}\left(  t\right)
:n_{0}\leq\ell\leq n\right\}  ,$ we obtain:
\[
R^{n+1}\left(  t\right)  \leq K+C\left\|  f_{0}\right\|  _{\infty}\int_{T}%
^{t}\left(  R^{n+1}\left(  s\right)  \right)  ^{3}ds+C\left(  t-T\right)
\;\;,\;\;t\geq T
\]
Let us select $\varepsilon_{0}=\frac{K}{C\left[  8\left\|  f_{0}\right\|
_{\infty}K^{3}+1\right]  }.$ It then follows, using a Gronwall type of
argument that:
\[
Q^{n}\left(  t\right)  \leq2K\;\;,\;\;n\geq n_{0}\;,\;\;0\leq t\leq
T+\varepsilon_{0}\;.
\]
\end{proof}

\bigskip\bigskip

\section{Energy estimate and consequences.}

The following energy estimates are standard for the Vlasov-Poisson system (cf.
\cite{Gl}). We state them here for further reference.

\begin{proposition}
Suppose that $f$ is a solution of (\ref{S1E1})-(\ref{S1E4a}) on the time
interval $t\in\left[  0,T\right]  $. Then:
\begin{equation}
\frac{d}{dt}\left(  \int_{\Omega\times\mathbb{R}^{3}}v^{2}fdxdv+\int_{\Omega
}\left(  E\right)  ^{2}dx\right)  =0\label{Eniter}%
\end{equation}
\end{proposition}

\begin{proof}
It is similar to the proof of the analogous result in the whole space (see
\cite{Gl}, page 120). The only difference is that we need to take into account
the contribution of some boundary terms. More precisely, using (\ref{S1E1}) we
obtain, after some integration by parts:
\begin{align}
\frac{d}{dt}\left(  \int_{\Omega\times\mathbb{R}^{3}}v^{2}fdxdv\right)   &
=2\int_{\Omega\times\mathbb{R}^{3}}\nabla_{x}\phi\cdot vfdxdv-2\int
_{\partial\Omega\times\mathbb{R}^{3}}v^{2}v\cdot n_{x}fdS_{x}dv\label{der1}\\
\frac{d}{dt}\left(  \int_{\Omega\times\mathbb{R}^{3}}E^{2}dx\right)   &
=-2\int_{\Omega\times\mathbb{R}^{3}}\nabla_{x}\phi\cdot vfdxdv+2\int
_{\partial\Omega}\phi\left(  \frac{\partial\phi}{\partial n_{x}}\right)
_{t}dS_{x}\label{der2}\\
&  +2\int_{\partial\Omega\times\mathbb{R}^{3}}\phi fv\cdot n_{x}%
dS_{x}dv\nonumber
\end{align}

The boundary term in (\ref{der1}) vanishes, since $\int_{\partial\Omega
\times\mathbb{R}^{3}}v^{2}v\cdot n_{x}fdS_{x}dv=0$ for each $x\in
\partial\Omega,$ due to the specular boundary condition (\ref{S1E4a}). On the
other hand, the two last terms in (\ref{der2}) vanish due to the fact that
$\frac{\partial\phi}{\partial n_{x}}=h$ is independent of $t$ and
$\int_{\mathbb{R}^{3}}fv\cdot n_{x}dx=0$ for each $x\in\partial\Omega.$ Adding
then (\ref{der1}), (\ref{der2}) we obtain (\ref{Eniter}).
\end{proof}

\bigskip

\begin{proposition}
Suppose that $f$ is a solution of (\ref{S1E1})-(\ref{S1E5new}) defined in
$0\leq t\leq T$ with $f\left(  0,x,v\right)  =f_{0}\left(  x,v\right)  $,
where $f_{0}$ is as in (\ref{S1E4a}). There exists $C$ depending only on $T$
and on the regularity norms assumed for $f_{0}$ in Theorem
\ref{globalexistence} such that:
\begin{align}
\sup_{0\leq t\leq T}\int_{\Omega}v^{2}f\left(  x,t\right)  dvdx  &  \leq
C\label{S2E1a}\\
\sup_{0\leq t\leq T}\left[  \left\|  \rho\left(  t,\cdot\right)  \right\|
_{L^{\frac{5}{3}}\left(  \Omega\right)  }\right]   &  \leq C\label{S2E1b}%
\end{align}%
\begin{equation}
\left\|  f\left(  t\right)  \right\|  _{L^{p}\left(  \Omega\times
\mathbb{R}^{3}\right)  }=\left\|  f_{0}\right\|  _{L^{p}\left(  \Omega
\times\mathbb{R}^{3}\right)  }\;,\text{ \ for all }1\leq p\leq\infty
.\label{S2E2b}%
\end{equation}

\begin{proof}
This result is just a consequence of (\ref{Eniter}) and its proof is standard
in kinetic theory. See for instance \cite{Gl}. Finally (\ref{S2E2b}) follows
multiplying (\ref{S1E1}) by $\left(  f\right)  ^{n-1},$ and integrating by
parts using the specular boundary conditions.
\end{proof}
\end{proposition}

\section{Pfaffelmoser's argument: Global bound for $Q\left(  t\right)  .$}

In this section we show that the function $Q\left(  t\right)  $ can be bounded
in any time interval $0\leq t\leq T$ and therefore that the corresponding
solutions of (\ref{S1E1})-(\ref{S1E5new}) can be extended to arbitrarily long
intervals. The estimate of $Q\left(  t\right)  $ will be derived using the
ideas of Pfaffelmoser (cf. \cite{Pfaffelmoser}) in the case of bounded domains
with purely reflected boundary conditions at $\partial\Omega.$\bigskip\ The
main content of the result is a uniform estimate for $Q\left(  t\right)  $ as
long as $f$ is defined.

Arguing as in the derivation of (\ref{rhonest}) we obtain the following
estimate:
\begin{equation}
\left\Vert \rho\right\Vert _{\infty}\leq\left\Vert f\right\Vert _{\infty
}Q\left(  t\right)  ^{3}\label{rhoestimate}%
\end{equation}
where $\rho$ is in (\ref{S1E2}).

The main result of this paper, that we will prove using Pfaffelmoser method,
is as follows:

\begin{theorem}
\label{boundforQ}Let $f_{0}\in C^{1,\mu}\left(  \Omega\times\mathbb{R}%
^{3}\right)  $ with $0<\mu<1.$ Suppose that $f\in C_{t,\left(  x,v\right)
}^{1;1,\lambda}\left(  [0,T]\times\Omega\times\mathbb{R}^{3}\right)  $ is a
solution of (\ref{S1E1})-(\ref{S1E4a}) with $\lambda\in\left(  0,1\right)
,\;0<T<\infty.$ There exists $\sigma\left(  T\right)  <\infty$ depending only
on $T,\;Q\left(  0\right)  $ and $\left\|  f_{0}\right\|  _{C^{1,\mu}\left(
\Omega\times\mathbb{R}^{3}\right)  }$ such that:
\begin{equation}
Q\left(  t\right)  \leq\sigma\left(  T\right)  \;\;,\;\;0\leq t\leq
T\;.\;\label{Qbound}%
\end{equation}
\end{theorem}

\subsection{Bounds for $Q\left(  t\right)  $}

Suppose that $\left(  \hat{X}\left(  s\right)  ,\hat{V}\left(  s\right)
\right)  $ is any fixed characteristic curve such that:
\[
\left(  \hat{X}\left(  0\right)  ,\hat{V}\left(  0\right)  \right)
\in\text{supp }f_{0}.
\]

The basic idea in the Pfaffelmoser's method consists in deriving estimates for
$E=\nabla\phi.$ Using the estimates in Proposition \ref{Neumann} it follows that:%

\begin{align}
\int_{t-\Delta}^{t}ds\left\vert E\left(  s,\hat{X}\left(  s\right)  \right)
\right\vert  &  \leq C\int_{t-\Delta}^{t}ds\int_{\Omega\times\mathbb{R}^{3}%
}\frac{f\left(  s,y,w\right)  }{\left\vert y-\hat{X}\left(  s\right)
\right\vert ^{2}}dydw+C\left\Vert h\right\Vert _{\infty}\Delta\label{SinE}\\
&  =C\int_{t-\Delta}^{t}ds\int_{\Omega\times\mathbb{R}^{3}}\frac{f\left(
t,x,v\right)  }{\left\vert X\left(  s\right)  -\hat{X}\left(  s\right)
\right\vert ^{2}}dxdv+C\left\Vert h\right\Vert _{\infty}\Delta\nonumber
\end{align}
where we are using the following change of variables:
\[
\left(  X\left(  s\right)  ,V\left(  s\right)  \right)  =\left(  y,w\right)
\rightarrow\left(  x,v\right)  =\left(  X\left(  t\right)  ,V\left(  t\right)
\right)
\]

We also note the measure preserving property ($dydw=dxdv$) that is due to the
fact that the evolution of the characteristics is hamiltonian away from the
boundary and that the measure $dxdv$ is also preserved by reflection on the boundary.

Pfaffelmoser's method is based on the idea of splitting the region of
integration in (\ref{SinE}) into three different sets that are usually termed
\textit{the good}, \textit{the bad,} and \textit{the ugly}. Fix $Q>0,$ that
will be precised later. In the rest of the argument we will use two numbers
$\Delta$ and $P$ defined by
\begin{equation}
P\equiv Q^{3/4-\delta}\;\;\;,\;\;\;\delta>0\;\;\text{small}\label{W3}%
\end{equation}%
\begin{equation}
\Delta\equiv\frac{c_{0}P}{Q^{2}}\label{W4}%
\end{equation}
where $c_{0}$ is small, but fixed number (independent on $Q$).

Given $Q,\;\Delta,\;P$ we define \textit{the good}, \textit{the bad,} and
\textit{the ugly} respectively by
\begin{align*}
G  &  \equiv\left\{  \left(  s,y,w\right)  \in\left[  t-\Delta,t\right]
\times\Omega\times\mathbb{R}^{3}:\;\left|  w\right|  \leq P,\text{ }\left|
w-\hat{V}\left(  t\right)  \right|  \leq P,\;\left|  w-\hat{V}^{+}\left(
t\right)  \right|  \leq P\right\} \\
B  &  \equiv\left\{  \left(  s,y,w\right)  \in\left[  t-\Delta,t\right]
\times\Omega\times\mathbb{R}^{3}:\left|  y-X\left(  s\right)  \right|
\leq\varepsilon_{0},\;\left|  w\right|  \geq P,\right. \\
&  \left.  \text{ }\left|  w-\hat{V}\left(  t\right)  \right|  \geq
P,\;\left|  w-\hat{V}^{+}\left(  t\right)  \right|  \geq P\right\} \\
U  &  \equiv\left\{  \left(  s,y,w\right)  \in\left[  t-\Delta,t\right]
\times\Omega\times\mathbb{R}^{3}:\left|  y-X\left(  s\right)  \right|
\geq\varepsilon_{0},\;\left|  w\right|  \geq P,\right.  \text{ }\\
&  \left.  \left|  w-\hat{V}\left(  t\right)  \right|  \geq P,\;\left|
w-\hat{V}^{+}\left(  t\right)  \right|  \geq P\right\}
\end{align*}
where
\[
\varepsilon_{0}\equiv\frac{R}{\left|  v\right|  ^{2}}\frac{1}{\left|
v-\hat{V}\left(  t\right)  \right|  }%
\]
if the characteristic curve $\hat{X}\left(  s\right)  $ does not intersect
$\partial\Omega$ on the interval $s\in\left[  t-\Delta,t\right]  ,$ and
otherwise
\[
\varepsilon_{0}\equiv\frac{R}{\left|  v\right|  ^{2}}\left(  \frac{1}{\left|
v-\hat{V}\left(  t\right)  \right|  }+\frac{1}{\left|  v-\hat{V}^{+}\left(
t\right)  \right|  }\right)  .
\]
Here
\[
\hat{V}^{+}\left(  t\right)  =\hat{V}^{\ast}\left(  s_{0}\right)  ,
\]
where
\[
s_{0}=\sup\left\{  s\in\left[  t-\Delta,t\right]  :\hat{X}\left(  s\right)
\in\partial\Omega\right\}  .
\]

The change due to the field can be estimated as:

\begin{lemma}
\label{Q4/3}Under the assumptions in Theorem \ref{boundforQ} we have the
following estimate:
\begin{equation}
\int_{t-\Delta}^{t}\left|  E\left(  X\left(  s\right)  ,s\right)  \right|
ds\leq C\Delta\left[  \left(  Q\left(  t\right)  \right)  ^{4/3}+1\right]
\;,\;\Delta\leq t\leq T\label{rough}%
\end{equation}
\end{lemma}

\begin{proof}
It is just an adaptation of a similar estimate in \cite{Gl}, page 122. \ The
last constant term in (\ref{rough}) is a consequence of the Neumann boundary
condition, i.e., $E=h$ at the boundary (cf. (\ref{S1E3})).
\end{proof}

\bigskip

We now estimate the change of $\hat{V}\left(  s\right)  $ due to geometry of
the domain and a basic geometric property that relates the change of $\hat
{V}\left(  s\right)  $ with the curvature will be investigated: The idea is as follows:

(i) The total length that the characteristic moves, including reflections is
bounded as $CQ\Delta.$

(ii) The change of the normal vector in a distance of order $\ell_{i}$ is
bounded by $C\ell_{i}.$

(iii) Therefore, the change of the angle of the vectors $\hat{V}\left(
s\right)  ,$ $\hat{V}^{\ast}\left(  s\right)  $ with respect to the previous
reflection is bounded by $\left(  CQ\ell_{i}+\int_{t_{i}}^{t_{i+1}}\left\vert
E\right\vert ds\right)  .$ (The change is the difference of outcoming vector
with respect to the next outcoming one, and incoming vectors with respect to
the incoming ones).

(iv) The total change of these vectors is then bounded by $CQ^{2}\Delta
+\int_{t-\Delta}^{t}\left\vert E\right\vert ds\leq C\left(  Q^{2}%
+Q^{4/3}\right)  \Delta\leq CQ^{2}\Delta,$ where the constant $C$ depends only
on the maximum of the curvature.

\bigskip

We now give its explicit formulation and the rigorous proof. To this end we
first recall basic ingredients such as Formulas of Frenet in differential
geometry (see \cite{S} p. 18, p. 94 for reference) that we will combine with
Lemma \ref{Festimate} in this paper.

\begin{lemma}
\label{DiffG}Let $T$ and $N$ be the unit tangential and normal vector on
$\partial\Omega$ respectively. Let $\kappa$ be the curvature along a curve and
$\kappa_{N}$ be the normal curvature in the direction $dx$ of the line of
curvature. Then we have
\[
dN+\kappa_{N}dx=0,\ dT=\kappa Nds,
\]

where $s$ is the arc length.
\end{lemma}

\begin{proof}
See \cite{S} p. 18, p. 94.
\end{proof}

In view of Lemma \ref{Festimate} the equations describing the evolution of the
characteristics in a geometrical form are:
\begin{align}
\frac{d\mu_{i}}{dt}  &  =\frac{w_{i}}{1+k_{i}x_{\perp}}\label{mu}\\
\frac{dw_{i}}{dt}  &  =E_{i}-\frac{v_{\perp}w_{i}k_{i}}{1+k_{i}x_{\perp}}%
-\sum_{j,\ell=1}^{2}\frac{\Gamma_{j,\ell}^{i}w_{j}w_{\ell}}{1+k_{j}x_{\perp}%
}\label{w}\\
\frac{dx_{\perp}}{dt}  &  =v_{\perp}\label{x-perp}\\
\frac{dv_{\perp}}{dt}  &  =F=E_{\perp}+\sum_{j=1}^{2}\frac{w_{j}^{2}b_{j}%
}{1+k_{j}x_{\perp}}\label{v-perp}%
\end{align}

Note that these equations work only near the boundary $\partial\Omega$, but
inside the domain the maximum displacement of the characteristics is $Q\Delta$
that can be made small.

\begin{lemma}
\label{VelGap}Let $\left(  X\left(  s\right)  ,V\left(  s\right)  \right)  $
and $\left(  \hat{X}\left(  s\right)  ,\hat{V}\left(  s\right)  \right)  $ be
two characteristics and let $Q$ be as in (\ref{Q}). Then we have
\[
\min\left\{  \left|  \hat{V}\left(  s\right)  -\hat{V}\left(  t\right)
\right|  ,\left|  \hat{V}\left(  s\right)  -\hat{V}^{+}\left(  t\right)
\right|  \right\}  \leq C\Delta Q^{2}\left(  t\right)  \;\;\;,\;\;\;s\in
\left[  t-\Delta,t\right]
\]%
\[
\min\left\{  \left|  V\left(  s\right)  -V\left(  t\right)  \right|  ,\left|
V\left(  s\right)  -V^{+}\left(  t\right)  \right|  \right\}  \leq C\Delta
Q^{2}\left(  t\right)  \;\;\;,\;\;\;s\in\left[  t-\Delta,t\right]
\]
\end{lemma}

\begin{proof}
Estimates equally apply to $V$ and $\hat{V}$ and so we only give a proof for
$V.$ Unlike the whole space case, we have to take care of the possible sign
change of $v_{\perp}$ at the bounces , i.e. if it becomes zero coming from the
region $v_{\perp}<0.$Using the equations (\ref{w}) and (\ref{v-perp}), it
follows that
\[
\left|  \frac{d\left(  \left|  v_{\perp}\right|  \right)  }{ds}\right|
+\left|  \frac{dw_{i}}{dt}\right|  \leq CQ^{2}\;,\;\;s\in\left[
t-\Delta,t\right]
\]
and this implies
\begin{equation}
\left|  w_{i}\left(  s\right)  -w_{i}\left(  t\right)  \right|  \leq
CQ^{2}\left(  t-s\right)  ,\label{wQ2}%
\end{equation}%
\begin{equation}
\left|  \left|  v_{\perp}\left(  s\right)  \right|  -\left|  v_{\perp}\left(
t\right)  \right|  \right|  \leq CQ^{2}\left(  t-s\right)  .\label{v-perpQ2}%
\end{equation}
(Notice that the equation above is valid even if we cross the bounces, and
that $\left|  v_{\perp}\right|  $ is differentiable). If the number of jumps
of the normal velocities is even, then by (\ref{wQ2})-(\ref{v-perpQ2}) and by
Lemma \ref{DiffG}, we have
\begin{align*}
\left|  V\left(  s\right)  -V\left(  t\right)  \right|   &  =\left|  \left(
\sum_{i=1}^{2}w_{i}u_{i}+v_{\perp}N\right)  \left(  s\right)  -\left(
\sum_{i=1}^{2}w_{i}u_{i}+v_{\perp}N\right)  \left(  t\right)  \right| \\
&  \leq\left|  \sum_{i=1}^{2}\left(  w_{i}\left(  s\right)  -w_{i}\left(
t\right)  \right)  u_{i}\left(  t\right)  +\left(  v_{\perp}\left(  s\right)
-v_{\perp}\left(  t\right)  \right)  N\left(  t\right)  \right| \\
&  +\sum_{i=1}^{2}\left|  w_{i}\left(  s\right)  \right|  \left|  u_{i}\left(
t\right)  -u_{i}\left(  s\right)  \right|  +\left|  v_{\perp}\left(  s\right)
\right|  \left|  N\left(  t\right)  -N\left(  s\right)  \right| \\
&  \leq CQ^{2}\left(  t-s\right)  .
\end{align*}
In a similar manner, we deduce that if the number of jumps of the normal
velocities is odd, then
\[
\left|  V\left(  s\right)  -V^{+}\left(  t\right)  \right|  \leq CQ^{2}\left(
t-s\right)  .
\]
Thus we conclude the assertion of the lemma. In general, we can split the time
interval $\left[  s,t\right]  $ into several time sub-intervals in such a way
that particles are governed completely by the equation (\ref{S2geom}) near the
boundary or by the usual Vlasov equation (\ref{S1E1}) away from the boundary
on each sub-interval. Then the estimates above combined with those in the
whole space and noticing that the effect $Q\left(  t\right)  ^{4/3}$ of the
electric field alone (Lemma \ref{Q4/3}) is negligible to that $Q\left(
t\right)  ^{2}$ of the geometry yield the lemma.
\end{proof}

\bigskip

\begin{lemma}
\label{BU}In the sets $B$ and $U,$ we have
\[
\frac{\left|  w\right|  }{2}\leq\left|  v\right|  \leq2\left|  w\right|  ,\
\]%
\[
\frac{\left|  w-\hat{V}\left(  t\right)  \right|  }{2}+\frac{\left|  w-\hat
{V}^{+}\left(  t\right)  \right|  }{2}\leq\left|  v-\hat{V}\left(  t\right)
\right|  +\left|  v-\hat{V}^{+}\left(  t\right)  \right|  \leq2\left|
w-\hat{V}\left(  t\right)  \right|  +2\left|  w-\hat{V}^{+}\left(  t\right)
\right|  .
\]
\end{lemma}

\begin{proof}
In the sets $B$ and $U$, we have $\left|  w\right|  \geq P,\left|  w-\hat
{V}\left(  t\right)  \right|  \geq P,\left|  w-\hat{V}^{+}\left(  t\right)
\right|  \geq P.$ By Lemma \ref{VelGap}, we have either $\left|  w-V\left(
t\right)  \right|  \leq C\Delta Q^{2}$ or $\left|  w-V^{+}\left(  t\right)
\right|  \leq C\Delta Q^{2}.$ If $\left|  w-V\left(  t\right)  \right|  \leq
C\Delta Q^{2},$ we have
\[
\left|  w-\hat{V}\left(  t\right)  \right|  -\left|  V\left(  t\right)
-w\right|  \leq\left|  v-\hat{V}\left(  t\right)  \right|  \leq\left|
w-\hat{V}\left(  t\right)  \right|  +\left|  V\left(  t\right)  -w\right|  ,
\]%
\[
\left|  w-\hat{V}^{+}\left(  t\right)  \right|  -\left|  V\left(  t\right)
-w\right|  \leq\left|  v-\hat{V}^{+}\left(  t\right)  \right|  \leq\left|
w-\hat{V}^{+}\left(  t\right)  \right|  +\left|  V\left(  t\right)  -w\right|
\]
Since, for small $c_{0},$%
\[
CQ^{2}\Delta\leq\frac{P}{4}\leq\frac{\left|  w-\hat{V}\left(  t\right)
\right|  }{4},
\]%
\[
CQ^{2}\Delta\leq\frac{P}{4}\leq\frac{\left|  w-\hat{V}^{+}\left(  t\right)
\right|  }{4},
\]
we have
\[
\frac{\left|  w-\hat{V}\left(  t\right)  \right|  }{2}\leq\left|  v-\hat
{V}\left(  t\right)  \right|  \leq2\left|  w-\hat{V}\left(  t\right)  \right|
,
\]%
\[
\frac{\left|  w-\hat{V}^{+}\left(  t\right)  \right|  }{2}\leq\left|
v-\hat{V}^{\ast}\left(  t\right)  \right|  \leq2\left|  w-\hat{V}^{+}\left(
t\right)  \right|  .
\]
In the other case, we similarly obtain
\[
\frac{\left|  w-\hat{V}^{+}\left(  t\right)  \right|  }{2}\leq\left|
v-\hat{V}\left(  t\right)  \right|  \leq2\left|  w-\hat{V}^{+}\left(
t\right)  \right|  ,
\]%
\[
\frac{\left|  w-\hat{V}\left(  t\right)  \right|  }{2}\leq\left|  v-\hat
{V}^{\ast}\left(  t\right)  \right|  \leq2\left|  w-\hat{V}\left(  t\right)
\right|  .
\]
Therefore, we deduce our lemma.
\end{proof}

\bigskip

We denote $X_{||}=\mu_{1}u_{1}+\mu_{2}u_{2},$ $X_{\perp}=x_{\perp}$ and
$V_{||}=w_{1}u_{1}+w_{2}u_{2},$ $V_{\perp}=v_{\perp}.$

\begin{lemma}
\label{more}If a trajectory $\left(  X,V\right)  $ has more than one bounce in
the interval $\left[  t-\Delta,t\right]  ,$ then we have, for all $s\in\left[
t-\Delta,t\right]  ,$%
\[
\left|  V_{\perp}\left(  s\right)  \right|  \leq CQ^{2}\left(  t\right)
\Delta.
\]
\end{lemma}

\begin{proof}
If a trajectory $\left(  X,V\right)  $ has more than one bounce, then we have
$V_{\perp}\left(  \tilde{s}\right)  =0,$ for some $\tilde{s}\in\left[
t-\Delta,t\right]  .$ Since
\begin{equation}
\left|  \frac{d\left|  V_{\perp}\right|  }{ds}\right|  \leq CQ^{2}\left(
t\right)  ,\label{normalVel}%
\end{equation}
the lemma follows.
\end{proof}

\bigskip

\begin{lemma}
\label{Xperp}Let $\left(  X\left(  s\right)  ,V\left(  s\right)  \right)  $
and $\left(  \hat{X}\left(  s\right)  ,\hat{V}\left(  s\right)  \right)  $ be
trajectories over $\left[  t-\Delta,t\right]  $. Suppose that
\[
\left|  X_{\perp}\left(  s_{0}\right)  -\hat{X}_{\perp}\left(  s_{0}\right)
\right|  =\min_{s\in\left[  t-\Delta,t\right]  }\left|  X_{\perp}\left(
s\right)  -\hat{X}_{\perp}\left(  s\right)  \right|  ,\ \ s_{0}\in\left(
t-\Delta,t\right)  .
\]
Then either both $X_{\perp}\left(  s_{0}\right)  >0,~\hat{X}_{\perp}\left(
s_{0}\right)  >0$ or both $X_{\perp}\left(  s_{0}\right)  =\hat{X}_{\perp
}\left(  s_{0}\right)  =0$.
\end{lemma}

\begin{proof}
We prove the lemma by contradiction. Suppose $X_{\perp}\left(  s_{0}\right)
>0 $ and $\hat{X}_{\perp}\left(  s_{0}\right)  =0.$ The case $\hat{X}_{\perp
}\left(  s_{0}\right)  >0$ and $X_{\perp}\left(  s_{0}\right)  =0$ can be
studied in a symmetric way. Notice that the function $\lambda\left(  s\right)
=\left|  X_{\perp}\left(  s\right)  -\hat{X}_{\perp}\left(  s\right)  \right|
^{2} $ is differentiable in the interval $\left(  t-\Delta,t\right)  $ except
at a finite set of points. Therefore, at the point $s=s_{0}$ where the minimum
of $\lambda$ is achieved we have:
\[
\frac{d\lambda}{ds}\left(  s_{0}-\right)  =\frac{d}{ds}\left|  X_{\perp
}\left(  s_{0}-\right)  -\hat{X}_{\perp}\left(  s_{0}-\right)  \right|
^{2}\leq0\text{ \ }%
\]
where from now on $f\left(  s_{0}-\right)  =\lim_{s\rightarrow s_{0},s<s_{0}%
}f\left(  s\right)  ,$ $f\left(  s_{0}+\right)  =\lim_{s\rightarrow
s_{0},s>s_{0}}f\left(  s\right)  .\;$We then have:
\[
X_{\perp}\left(  s_{0}\right)  \left(  V_{\perp}\left(  s_{0}-\right)
-\hat{V}_{\perp}\left(  s_{0}-\right)  \right)  \leq0,
\]
which implies
\begin{equation}
V_{\perp}\left(  s_{0}-\right)  \leq\hat{V}_{\perp}\left(  s_{0}-\right)
<0.\label{neg}%
\end{equation}

The fact that $\hat{V}_{\perp}\left(  s_{0}-\right)  \neq0$ is a consequence
of Lemma \ref{alphbound}.

Notice that, since $\hat{X}_{\perp}\left(  s_{0}\right)  =0$ there is a
reflection of $\hat{V}_{\perp}$ at $s=s_{0}.$ On the other hand, since
$X_{\perp}\left(  s_{0}\right)  >0$ ,\ $V_{\perp}$ is continuous at $s=s_{0}%
.$\ Therefore:
\[
V_{\perp}\left(  s_{0}+\right)  =V_{\perp}\left(  s_{0}-\right)
\;\;\;,\;\;\hat{V}_{\perp}\left(  s_{0}+\right)  =-\hat{V}_{\perp}\left(
s_{0}-\right)
\]

Thus we have
\[
\left|  X_{\perp}\left(  s\right)  -\hat{X}_{\perp}\left(  s\right)  \right|
=X_{\perp}\left(  s\right)  -\hat{X}_{\perp}\left(  s\right)  =X_{\perp
}\left(  s_{0}\right)  +\left(  V_{\perp}\left(  s_{0}-\right)  +\hat
{V}_{\perp}\left(  s_{0}-\right)  \right)  \left(  s-s_{0}\right)  +o\left(
s-s_{0}\right)
\]
as $s\rightarrow s_{0},\;s>s_{0}.\;$Due to (\ref{neg}), we have $\left|
X_{\perp}\left(  s\right)  -\hat{X}_{\perp}\left(  s\right)  \right|
<X_{\perp}\left(  s_{0}\right)  =\left|  X_{\perp}\left(  s_{0}\right)
-\hat{X}_{\perp}\left(  s_{0}\right)  \right|  $ for $s-s_{0}>0$ sufficiently
small, but this contradicts the fact that $\left|  X_{\perp}\left(  s\right)
-\hat{X}_{\perp}\left(  s\right)  \right|  $ reaches its minimum at $s=s_{0}$.
Therefore the result follows.
\end{proof}

We show the following crucial separation property.

\begin{lemma}
(Separation property) In the ugly set $U$, there exist $s_{0},s_{1}\in\left[
t-\Delta,t\right]  $ such that the following separation holds:
\begin{equation}
\left|  X\left(  s\right)  -\hat{X}\left(  s\right)  \right|  \geq C\left(
\varepsilon_{0}+\min\left\{  \left|  v-\hat{V}\left(  t\right)  \right|
\left|  s-s_{0}\right|  ,\left|  v-\hat{V}^{+}\left(  t\right)  \right|
\left|  s-s_{1}\right|  \right\}  \right)  ,\;\;s\in\left[  t-\Delta,t\right]
\label{S2E3}%
\end{equation}
where $C$ is a universal constant depending only on the curvature of
$\partial\Omega.$
\end{lemma}

\begin{proof}
We separate into two cases:

\textbf{Case 1}: \ Both trajectories $\left(  X\left(  s\right)  ,V\left(
s\right)  \right)  ,$ $\left(  \hat{X}\left(  s\right)  ,\hat{V}\left(
s\right)  \right)  $ have at most one bounce in the time interval $\left[
t-\Delta,t\right]  $. Let $t-\Delta\leq t_{1}<t_{2}\leq t$ be possible two
bouncing times with $t_{1},t_{2}$ corresponding to $\left(  X\left(  s\right)
,V\left(  s\right)  \right)  ,$ $\left(  \hat{X}\left(  s\right)  ,\hat
{V}\left(  s\right)  \right)  $ respectively. We split the time interval
$\left[  t-\Delta,t\right]  $ into a maximum of three sub-intervals, namely
$\left[  t-\Delta,t_{1}\right]  \cup\left[  t_{1},t_{2}\right]  \cup\left[
t_{2},t\right]  .$ In the most general case, some of these intervals could be
empty. Let us describe the argument in the most general case, since for a
smaller number of reflections the argument required is just a minor
simplification of it. Pick $s_{0}$ and $s_{1}$ such that
\begin{align*}
\min_{s\in\left[  t-\Delta,t_{1}\right]  \cup\lbrack t_{2},t]}\left|  \left(
X-\hat{X}\right)  \left(  s\right)  \right|   &  =\left|  \left(  X-\hat
{X}\right)  \left(  s_{0}\right)  \right|  ,\\
\min_{s\in\lbrack t_{1},t_{2}]}\left|  \left(  X-\hat{X}\right)  \left(
s\right)  \right|   &  =\left|  \left(  X-\hat{X}\right)  \left(
s_{1}\right)  \right|  .
\end{align*}
In the interval $[t_{1},t_{2}]$ there is no bounces along the trajectories.
Then, we can argue exactly as in the case without boundaries (cf. \cite{Gl},
pages 128-129) to show that:
\[
\left|  X\left(  s\right)  -\hat{X}\left(  s\right)  \right|  \geq C\left|
v-\hat{V}^{+}\left(  t\right)  \right|  \left|  s-s_{1}\right|  ,\text{ \ for
}s\in\left[  t_{1},t_{2}\right]  .
\]

On the other hand, the portion of the trajectories $X\left(  s\right)
,\hat{X}\left(  s\right)  $ for $s\in\lbrack t_{2},t]$ might be reflected with
respect to the boundary $\partial\Omega.$ Suppose for the moment that the
boundary $\partial\Omega$ is flat as in the half space case. The trajectories
over $[t_{2},t]$ obtained by reflections together with the original
trajectories $X\left(  s\right)  ,\;\hat{X}\left(  s\right)  $ for
$s\in\left[  t-\Delta,t_{1}\right]  $ yield portion of new trajectories
without bounces, and satisfying an equation of the form:
\begin{align*}
\frac{dX}{ds}  &  =V\\
\frac{dV}{ds}  &  =\tilde{E}%
\end{align*}
where $\int_{t-\Delta}^{t}\left|  \tilde{E}\left(  X\left(  s\right)
,s\right)  \right|  ds\leq C\Delta\left[  \left(  Q\left(  t\right)  \right)
^{4/3}+1\right]  .$ We can argue then exactly as in the case of the whole
space (cf. \cite{Gl}), to estimate the difference between the trajectories,
and since the reflection with respect to the plane $x_{1}=0$ is an isometry,
we finally obtain:
\begin{equation}
\left|  X\left(  s\right)  -\hat{X}\left(  s\right)  \right|  \geq C\left|
V\left(  t\right)  -\hat{V}\left(  t\right)  \right|  \left|  s-s_{0}\right|
,\text{ \ for }s\in\left[  t-\Delta,t_{1}\right]  \cup\left[  t_{2},t\right]
.\label{sepCase1}%
\end{equation}

Now if the boundary $\partial\Omega$ is not flat then the change of the normal
vectors between two reflections can be bounded by $C\Delta Q$ and the
corresponding change of the vectors $V$ and $\hat{V}$ can be bounded by
$C\Delta Q^{2},$ which is smaller than $c_{0}P$ for sufficiently small $c_{0}.
$ Therefore the change of these vectors is small compared to $\left|  V\left(
t\right)  -\hat{V}\left(  t\right)  \right|  $ in the ugly set and the
inequality (\ref{sepCase1}) above is valid for the case of domain with curvature.

\textbf{Case 2}: At least one of the trajectories has more than one bounce in
$\left[  t-\Delta,t\right]  .$ Let $\left(  \hat{X},\hat{V}\right)  $ have
more than one bounce in $\left[  t-\Delta,t\right]  .$

We first consider the case
\begin{equation}
\left|  w_{i}\left(  s\right)  -\hat{w}_{i}\left(  s\right)  \right|
\geq\frac{1}{2}\left|  V\left(  s\right)  -\hat{V}\left(  s\right)  \right|
,\ \text{for all }s\in\left[  t-\Delta,t\right]  ,\label{Tan}%
\end{equation}

i.e., the tangential part of $V-\hat{V}$ dominates along the trajectory in the
whole interval $\left[  t-\Delta,t\right]  .$ Let
\[
\min_{s\in\lbrack t-\Delta,t]}\left|  \mu_{i}\left(  s\right)  -\hat{\mu}%
_{i}\left(  s\right)  \right|  =\left|  \mu_{i}\left(  s_{0}\right)  -\hat
{\mu}_{i}\left(  s_{0}\right)  \right|  .
\]
Note that the tangential part of the trajectory, $X_{||}$ is $C^{1}$. Consider
the equation (\ref{mu}) and (\ref{w}) for the tangential components of the
position and the velocity:
\[
\frac{dw_{i}}{dt}=E_{i}-\frac{v_{\perp}w_{i}k_{i}}{1+k_{i}x_{\perp}}%
-\sum_{j,\ell=1}^{2}\frac{\Gamma_{j,\ell}^{i}w_{j}w_{\ell}}{1+k_{j}x_{\perp}%
}=\mathcal{O}\left(  Q^{4/3}\right)  +\mathcal{O}\left(  Q^{2}\right)  ,
\]%
\[
\frac{d\mu_{i}}{ds}=\frac{w_{i}}{1+k_{i}x_{\perp}}=w_{i}\left(  s_{0}\right)
+\mathcal{O}\left(  \Delta Q^{2}\right)  +O\left(  x_{\perp}\left|
w_{i}\right|  \right)  =w_{i}\left(  s_{0}\right)  +\mathcal{O}\left(  \Delta
Q^{2}\right)  ,
\]
where the second term comes from the change on the velocity $w_{i}$ due to the
field and to the geometry and constants depend only on the geometry of the
domain. We also used the fact that
\[
\left|  x_{\perp}\right|  \left|  w_{i}\right|  =\mathcal{O}\left(  \Delta
Q^{2}\right)  .
\]

Then, integrating the tangential part $\mu_{i}-\hat{\mu}_{i}$ of the
difference of $\left(  X,V\right)  $ and $\left(  \hat{X},\hat{V}\right)  $,
we get:
\[
\left(  \mu_{i}-\hat{\mu}_{i}\right)  \left(  s\right)  =\left(  \mu_{i}%
-\hat{\mu}_{i}\right)  \left(  s_{0}\right)  +\left(  w_{i}\left(
s_{0}\right)  -\hat{w}_{i}\left(  s_{0}\right)  \right)  \left(
s-s_{0}\right)  +O\left(  \Delta Q^{2})\left(  s-s_{0}\right)  \right)
\]
By Lemma \ref{VelGap} and (\ref{Tan}), we have
\[
\left|  w_{i}\left(  s_{0}\right)  -\hat{w}_{i}\left(  s_{0}\right)  \right|
\geq\frac{1}{2}\left|  V\left(  s_{0}\right)  -\hat{V}\left(  s_{0}\right)
\right|  \geq\frac{1}{2}\min\left\{  \left|  v-\hat{V}\left(  t\right)
\right|  ,\left|  v-\hat{V}^{+}\left(  t\right)  \right|  \right\}
+\mathcal{O}\left(  Q^{2}\left(  t\right)  \Delta\right)  .
\]
Since $\left|  \left(  \mu_{i}-\hat{\mu}_{i}\right)  \left(  s\right)
\right|  ^{2}$ is a $C^{1}$ function it folllows that at the point $s_{0}%
\in\lbrack t-\Delta,t]$ where $\left|  \left(  \mu_{i}-\hat{\mu}_{i}\right)
\left(  s\right)  \right|  ^{2}$ attains the minimum we have
\begin{align*}
&  \left(  s-s_{0}\right)  \frac{d}{ds}\left(  \left|  \left(  \mu_{i}%
-\hat{\mu}_{i}\right)  \left(  s\right)  \right|  ^{2}\right) \\
&  =[\mu_{i}\left(  s_{0}\right)  -\hat{\mu}_{i}\left(  s_{0}\right)
]\cdot\lbrack w_{i}\left(  s_{0}\right)  -\hat{w}_{i}\left(  s_{0}\right)
]\left(  s-s_{0}\right)  \geq0.
\end{align*}
Note that the inequality above takes into account the possibility that the
minimun of $\left|  \left(  \mu_{i}-\hat{\mu}_{i}\right)  \left(  s\right)
\right|  ^{2}$ can be achieved at the end points $s_{0}=t-\Delta,t.$ Since
\[
\mathcal{O}\left(  Q^{2}\left(  t\right)  \Delta\right)  =c_{0}P,
\]
we deduce
\[
\left|  \left(  \mu_{i}-\hat{\mu}_{i}\right)  \left(  s\right)  \right|
\geq\frac{1}{4}\min\left\{  \left|  v-\hat{V}\left(  t\right)  \right|
,\left|  v-\hat{V}^{+}\left(  t\right)  \right|  \right\}  \left(
s-s_{0}\right)  ,
\]
for sufficiently small $c_{0}$. Using the definition of the ugly set we have
$\left|  \left(  \mu_{i}-\hat{\mu}_{i}\right)  \left(  s\right)  \right|
\geq\varepsilon_{0}$. Thus (\ref{S2E3}) follows.

We now consider the complementary case, i.e., there is $\bar{s}\in\left[
t-\Delta,t\right]  $ such that
\[
\left|  V_{\perp}\left(  \bar{s}\right)  -\hat{V}_{\perp}\left(  \bar
{s}\right)  \right|  \geq\frac{1}{2}\left|  V\left(  \bar{s}\right)  -\hat
{V}\left(  \bar{s}\right)  \right|  .
\]
By (\ref{normalVel}), we have
\[
\left|  \left|  V_{\perp}\left(  s\right)  \right|  -\left|  \hat{V}_{\perp
}\left(  s\right)  \right|  \right|  \geq\left|  \left|  V_{\perp}\left(
\bar{s}\right)  \right|  -\left|  \hat{V}_{\perp}\left(  \bar{s}\right)
\right|  \right|  -CQ^{2}\left(  t\right)  \Delta.
\]
On the other hand, note that
\[
\left|  \left|  V_{\perp}\left(  s\right)  \right|  -\left|  \hat{V}_{\perp
}\left(  s\right)  \right|  \right|  =\left|  V_{\perp}\left(  s\right)
-sgn\left(  V_{\perp}\left(  s\right)  \right)  \left|  \hat{V}_{\perp}\left(
s\right)  \right|  \right|  .
\]
Using Lemma \ref{more}, it follows that
\begin{align*}
&  \left|  V_{\perp}\left(  s\right)  -sgn\left(  V_{\perp}\left(  s\right)
\right)  \left|  \hat{V}_{\perp}\left(  s\right)  \right|  \right| \\
&  =\left|  V_{\perp}\left(  s\right)  -\hat{V}_{\perp}\left(  s\right)
+\hat{V}_{\perp}\left(  s\right)  -sgn\left(  V_{\perp}\left(  s\right)
\right)  \left|  \hat{V}_{\perp}\left(  s\right)  \right|  \right| \\
&  \leq\left|  V_{\perp}\left(  s\right)  -\hat{V}_{\perp}\left(  s\right)
\right|  +2\left|  \hat{V}_{\perp}\left(  s\right)  \right| \\
&  \leq\left|  V_{\perp}\left(  s\right)  -\hat{V}_{\perp}\left(  s\right)
\right|  +CQ^{2}\left(  t\right)  \Delta.
\end{align*}
Similarly, we get
\begin{align*}
&  \left|  V_{\perp}\left(  \bar{s}\right)  -sgn\left(  V_{\perp}\left(
\bar{s}\right)  \right)  \left|  \hat{V}_{\perp}\left(  \bar{s}\right)
\right|  \right| \\
&  =\left|  V_{\perp}\left(  \bar{s}\right)  -\hat{V}_{\perp}\left(  \bar
{s}\right)  +\hat{V}_{\perp}\left(  \bar{s}\right)  -sgn\left(  V_{\perp
}\left(  \bar{s}\right)  \right)  \left|  \hat{V}_{\perp}\left(  \bar
{s}\right)  \right|  \right| \\
&  \geq\left|  V_{\perp}\left(  \bar{s}\right)  -\hat{V}_{\perp}\left(
\bar{s}\right)  \right|  -2\left|  \hat{V}_{\perp}\left(  \bar{s}\right)
\right| \\
&  \geq\left|  V_{\perp}\left(  \bar{s}\right)  -\hat{V}_{\perp}\left(
\bar{s}\right)  \right|  -CQ^{2}\left(  t\right)  \Delta.
\end{align*}
Thus, we obtain, for all $s\in\left[  t-\Delta,t\right]  ,$%
\begin{align}
\left|  V_{\perp}\left(  s\right)  -\hat{V}_{\perp}\left(  s\right)  \right|
&  \geq\left|  V_{\perp}\left(  \bar{s}\right)  -\hat{V}_{\perp}\left(
\bar{s}\right)  \right|  -CQ^{2}\left(  t\right)  \Delta\label{lowineq}\\
&  \geq\frac{1}{2}\left|  V\left(  \bar{s}\right)  -\hat{V}\left(  \bar
{s}\right)  \right|  -CQ^{2}\left(  t\right)  \Delta\nonumber\\
&  \geq\frac{1}{2}\min\{\left|  v-\hat{V}\left(  t\right)  \right|  ,\left|
v-\hat{V}^{+}\left(  t\right)  \right|  \}-CQ^{2}\left(  t\right)
\Delta\nonumber\\
&  \geq\frac{1}{2}\min\{\left|  v-\hat{V}\left(  t\right)  \right|  ,\left|
v-\hat{V}^{+}\left(  t\right)  \right|  \}-c_{0}P\nonumber\\
&  \geq\frac{1}{4}\min\{\left|  v-\hat{V}\left(  t\right)  \right|  ,\left|
v-\hat{V}^{+}\left(  t\right)  \right|  \}\geq\frac{P}{4}.\nonumber
\end{align}
Using Lemma \ref{more}, it follows, choosing $c_{0}$ in (\ref{W4})
sufficiently small, that:
\begin{equation}
\left|  V_{\perp}\left(  s\right)  \right|  \geq\frac{P}{8}\label{lowest}%
\end{equation}
for all $s\in\left[  t-\Delta,t\right]  .$ Taking into account (\ref{rough})
it follows that $V_{\perp}\left(  s\right)  $ changes sign, by reflection, at
most once in the interval $s\in\left[  t-\Delta,t\right]  $ if $c_{0}$ is
sufficiently small. Combining Lemma \ref{more}, (\ref{lowineq}), and
(\ref{lowest}) it follows that $V_{\perp}\left(  s\right)  -\hat{V}_{\perp
}\left(  s\right)  $ changes sign at most once for $s\in\left[  t-\Delta
,t\right]  .$ Indeed, $\hat{V}_{\perp}\left(  s\right)  $ is small compared to
$\left|  V_{\perp}\left(  s\right)  \right|  \geq\frac{P}{8}$ in the interval
$\left[  t-\Delta,t\right]  ,$ and $V_{\perp}\left(  s\right)  $ changes sign
only once at most. Suppose that $V_{\perp}\left(  s\right)  $ changes sign at
$s=s_{0}.$ Since $X_{\perp}\left(  s\right)  \geq0,$ it follows that
$\operatorname*{sign}\left(  V_{\perp}\left(  s\right)  \right)
=\operatorname*{sign}\left(  V_{\perp}\left(  s\right)  -\hat{V}_{\perp
}\left(  s\right)  \right)  =\operatorname*{sign}\left(  s-s_{0}\right)
,\;$for $s\in\left[  t-\Delta,t\right]  .$ We have\
\[
X_{\perp}\left(  s\right)  -\hat{X}_{\perp}\left(  s\right)  =\int_{s_{0}}%
^{s}\left[  V_{\perp}\left(  \tau\right)  -\hat{V}_{\perp}\left(  \tau\right)
\right]  d\tau.
\]
Then, using (\ref{lowineq}), we have:
\[
\left|  X_{\perp}\left(  s\right)  -\hat{X}_{\perp}\left(  s\right)  \right|
\geq\frac{1}{4}\min\{\left|  V\left(  t\right)  -\hat{V}\left(  t\right)
\right|  ,\left|  V\left(  t\right)  -\hat{V}^{+}\left(  t\right)  \right|
\}\left|  s-s_{0}\right|  .
\]
Since in the ugly set $U,$ $\left|  X_{\perp}\left(  s\right)  -\hat{X}%
_{\perp}\left(  s\right)  \right|  \geq\varepsilon_{0},\;$we obtain
(\ref{S2E3}). The proof is complete.
\end{proof}

\bigskip

\begin{proof}
[\textbf{Proof of Theorem \ref{boundforQ}}]The key point in the proof is to
estimate the right-hand side of (\ref{SinE}). Let us assume without loss of
generality that $Q\geq1,$ since otherwise the leading contribution would be
$C\Delta$ in (\ref{SinE}). In order to make this estimate we separate the
contributions of the sets\ $G,\;B$ and $U$:
\begin{align*}
&  \int_{t-\Delta}^{t}ds\int_{\Omega\times IR^{3}}\frac{f\left(  y,w,s\right)
}{\left|  y-\hat{X}\left(  s\right)  \right|  ^{2}}dydw\\
&  =\int_{G}\left[  ...\right]  dsdydw+\int_{B}\left[  ...\right]
dsdydw+\int_{U}\left[  ...\right]  dsdydw
\end{align*}
In order to estimate the contribution of the good set we define:
\[
\rho_{G}\left(  y,s\right)  \equiv\int_{G}f\left(  y,w,s\right)  dw
\]
Standard estimates yield
\[
\left\|  \rho_{G}\right\|  _{\infty}\leq CP^{3}%
\]
where from now on $C$ depends on $\left\|  f_{0}\right\|  _{\infty},$ but not
on $P,\;\varepsilon_{0},\;R,\;Q$,$\;\Delta.$ Arguing as in the derivation of
(\ref{rough}) we obtain:
\begin{equation}
\int_{G}\left[  ...\right]  dsdydw\leq\Delta P^{4/3}\label{good}%
\end{equation}
In order to estimate the contribution of the bad set, notice that Lemma
\ref{BU} implies
\[
\varepsilon_{0}\leq\frac{8R}{\left|  w\right|  ^{2}}\left(  \frac{1}{\left|
w-\hat{V}\left(  t\right)  \right|  }+\frac{1}{\left|  w-\hat{V}^{+}\left(
t\right)  \right|  }\right)  .
\]
Then
\[
\int_{B}\frac{f\left(  y,w,s\right)  }{\left|  y-\hat{X}\left(  s\right)
\right|  ^{2}}dsdydw\leq C\int_{t-\Delta}^{t}ds\int_{\left|  w\right|  \leq
Q}\varepsilon_{0}dw=
\]%
\begin{equation}
C\int_{t-\Delta}^{t}ds\int_{\left|  w\right|  \leq Q}\frac{R}{\left|
w\right|  ^{2}}\left(  \frac{1}{\left|  w-\hat{V}\left(  t\right)  \right|
}+\frac{1}{\left|  w-\hat{V}^{+}\left(  t\right)  \right|  }\right)  dw\leq
CR\Delta\log\left(  Q\right) \label{bad}%
\end{equation}
where by assumption $\left|  \hat{V}\left(  t\right)  \right|  \geq1,$ since
otherwise the corresponding characteristic does not have effect in the
variation of $Q.$

Lastly we estimate the integral over the ugly set:

The aim is to estimate
\[
\int_{U}\frac{f\left(  y,w,s\right)  }{\left|  y-\hat{X}\left(  s\right)
\right|  ^{2}}dsdydw=\int_{U}\frac{f\left(  x,v,t\right)  }{\left|  X\left(
s\right)  -\hat{X}\left(  s\right)  \right|  ^{2}}dsdxdv
\]
Using (\ref{S2E3}) we can estimate this integral as
\[
C\int_{U}\frac{f\left(  x,v,t\right)  }{\left(  \varepsilon_{0}+\min\left\{
\left|  v-\hat{V}\left(  t\right)  \right|  \left|  s-s_{0}\right|  ,\left|
v-\hat{V}^{+}\left(  t\right)  \right|  \right\}  \left|  s-s_{1}\right|
\right)  ^{2}}dsdxdv
\]
The integration with respect to time $s$ can be estimated as
\begin{align*}
&  \int_{t-\Delta}^{t}\frac{ds}{\left(  \varepsilon_{0}+\min\left\{  \left|
v-\hat{V}\left(  t\right)  \right|  \left|  s-s_{0}\right|  ,\left|  v-\hat
{V}^{+}\left(  t\right)  \right|  \left|  s-s_{1}\right|  \right\}  \right)
^{2}}\\
&  \leq\frac{C}{\varepsilon_{0}}\left(  \frac{1}{\left|  v-\hat{V}\left(
t\right)  \right|  }+\frac{1}{\left|  v-\hat{V}^{\ast}\left(  t\right)
\right|  }\right) \\
&  \leq\frac{C}{\varepsilon_{0}}\left(  \frac{1}{\left|  v-\hat{V}\left(
t\right)  \right|  }+\frac{1}{\left|  v-\hat{V}^{\ast}\left(  t\right)
\right|  }\right)  =\frac{Cv^{2}}{R}%
\end{align*}
Then
\begin{equation}
\int_{U}\left[  ...\right]  dsdydw\leq\frac{C}{R}\int v^{2}f\left(
x,v,t\right)  dvdx\leq\frac{C}{R\Delta}\Delta\label{ugly}%
\end{equation}
Combining the estimates (\ref{good}), (\ref{bad}), and (\ref{ugly}) for the
Good, Bad, and Ugly set, we obtain
\begin{align*}
\int_{t-\Delta}^{t}\left|  E\left(  s,\hat{X}\left(  s\right)  \right)
\right|  ds  &  \leq C\Delta\left(  P^{4/3}+R\log\left(  Q\right)  +\frac
{1}{R\Delta}\right) \\
&  =C\Delta\left(  P^{4/3}+R\log\left(  Q\right)  +\frac{Q^{4/3}}{RP}\right)
\end{align*}
Choosing $R=Q^{1-\delta},$ $P=Q^{3/4-\delta}$ we obtain
\[
\int_{t-\Delta}^{t}\left|  E\left(  s,\hat{X}\left(  s\right)  \right)
\right|  ds\leq C\Delta Q^{\beta}%
\]
where $\beta<1.$

Therefore
\[
\frac{Q\left(  t\right)  -Q\left(  t-\Delta\right)  }{\Delta}\leq C\left(
Q\left(  t\right)  \right)  ^{\beta},
\]
and a standard\ iteration yields $Q\left(  t\right)  $ bounded in any interval
$0\leq t\leq T.$ Thus this completes the proof.
\end{proof}

\bigskip

\begin{proof}
[\textbf{Proof of Theorem \ref{globalexistence}}]In order to prove the
existence of a solution $f$ of (\ref{S1E1})-(\ref{S1E4a}) globally defined in
$t$ we will show that the sequence of functions $f^{n}$ defined by
(\ref{itf0})-(\ref{itfield}) converges as $n\rightarrow\infty$ to a solution
of (\ref{itfn1})-(\ref{itfn5}) for arbitrary values of $t.$ To this end, it
suffices to show that the functions $Q^{n}\left(  t\right)  $ are uniformlly
bounded in each compact set of $t\in\left[  0,\infty\right)  .$The desired
limit property would then follow from Proposition \ref{Propconv}.

To this end, we define $L\left(  t\right)  =\sup_{n}Q^{n}\left(  t\right)  .$
We have that $L\left(  \cdot\right)  $ is increasing in $t.$ We denote as
$T_{\max}$ the time where:
\[
\lim_{t\rightarrow T_{\max}}L\left(  t\right)  =\infty
\]
Our goal is to show that $T_{\max}=\infty.$ Let us assume that $T_{\max
}<\infty.$ We define $\varepsilon_{0}=\varepsilon_{0}\left(  2\sigma\left(
T_{\max}\right)  ,\left\|  f_{0}\right\|  _{\infty}\right)  ,$ where the
function $\varepsilon_{0}\left(  \cdot\right)  $ is as in Proposition
\ref{Prolog} and the function $\sigma\left(  T\right)  $ is as in Theorem
\ref{boundforQ}. Notice that the functions $Q^{n}\left(  t\right)  $ are
uniformly bounded for $t\in\left[  0,T_{\max}-\frac{\varepsilon_{0}}%
{2}\right]  $\ by definition of $T_{\max}.$ Therefore, Proposition
\ref{Prolog} implies that $f^{n}\rightarrow f$ in $C_{t,\left(  x,v\right)
}^{\nu;1,\lambda}$ for $0\leq t\leq T_{\max}-\frac{\varepsilon_{0}}{2}%
,\;0<\nu<1.$ We can use Proposition \ref{Qconv} to prove that $Q^{n}\left(
t\right)  \rightarrow Q\left(  t\right)  $ for $t\in\left[  0,T_{\max}%
-\frac{\varepsilon_{0}}{2}\right]  .$ In particular, $\lim_{n\rightarrow
\infty}Q_{n}\left(  \bar{t}\right)  =Q\left(  \bar{t}\right)  \leq
\sigma\left(  T_{\max}\right)  ,$ for $\bar{t}=T_{\max}-\frac{\varepsilon_{0}%
}{2}.$ Therefore $Q^{n}\left(  \bar{t}\right)  \leq2\sigma\left(  T_{\max
}\right)  $ for $n\geq n_{0}$ with $n_{0}$ large. Then, Proposition
\ref{Prolog}\ implies that the sequence $Q^{n}\left(  t\right)  $ is uniformly
bounded for $t\in\left[  0,T_{\max}+\frac{\varepsilon_{0}}{2}\right]  ,$
whence $L\left(  t\right)  $ is bounded as $t\rightarrow T_{\max}.$ This
contradicts the definition of $T_{\max}$ and concludes the proof of the
existence of a solution of (\ref{S1E1})-(\ref{S1E4a}) in $C_{t,\left(
x,v\right)  }^{1;1,\lambda}$ for some $0<\lambda<\mu$ and for $0\leq t\leq T$
as asserted in Theorem \ref{globalexistence}. In order to prove uniqueness we
argue as in the proof of (\ref{gronw1}) to obtain, that two $C_{t,\left(
x,v\right)  }^{1;1,\lambda}$ solutions of (\ref{S1E1})-(\ref{S1E4a}) with the
same initial and boundary satisfy
\[
\left\|  f_{1}\left(  t\right)  -f_{2}\left(  t\right)  \right\|  _{L^{1}}\leq
C\left(  T\right)  \int_{0}^{t}\left\|  f_{1}\left(  s\right)  -f_{2}\left(
s\right)  \right\|  _{L^{1}}ds,
\]
Therefore $f_{1}=f_{2}.$ This completes the proof of the Theorem.
\end{proof}

\bigskip

\begin{acknowledgement}
The authors thank the hospitality and the support of the Max Planck Institut
for Mathematics in the Sciences in Leipzig. HJH is supported by Korean
Research Fund KRF-2007-331-C00021, Postech Research Fund 2007 and PMI. JJLV
acknowledges also support of the Alexander von Humboldt foundation and DGES
Grant MTM2004-05634.
\end{acknowledgement}


\bigskip

\begin{thebibliography}{99}
\bibitem{BD}C. Bardos, D. Degond, Global existence for the Vlasov-Poisson
equation in 3 space variables with small initial data, Ann. Inst. H.
Poincar\'{e}. Analyse non lineaire 2(1985), pp. 101-118.

\bibitem {BR}J. Batt, G. Rein, Global classical solutions of the periodic
Vlasov-Poisson system in three dimensions, C. R. Acad. Sci. Paris, 313(1991),
pp. 411-416.

\bibitem {E}D. M. Eidus, Inequalities for Green's function, Mat. Sb.,
87(1958), pp. 455-470 [Russian].

\bibitem {GT}D. Gilbarg, N. S. Trudinger, \textit{Elliptic Partial
Differential Equations of Second Order.} Springer, Berlin, 1983.

\bibitem {Gl}R. Glassey, The Cauchy problem in kinetic theory, SIAM:
Philadelphia, PA, 1996.

\bibitem {YanGuo}Y. Guo, Singular solutions of Vlasov-Maxwell boundary
problems in one dimension. Arch. Rational Mech. Anal. \textbf{131}, 241-304 (1995).

\bibitem {GuoIndiana}Y. Guo, Regularity for the Vlasov equations in a Half
space. Indiana Univ. Math. J. \textbf{43}, 255-320 (1994).

\bibitem {GS1}R. Glassey, W. Strauss, Singularity formation in a collisionless
plasma could occur only at high velocities, Arch. Rational Mech. Anal.,
92(1986), pp. 59-90.

\bibitem {Hwang}H.J. Hwang, Regularity for the Vlasov-Poisson system in a
convex domain. SIAM J. Math. Anal. \textbf{36} , 121-171 (2004).

\bibitem {HV}H.J. Hwang and J. J. L. Vel\'{a}zquez, A new proof \ of global
existence for the Vlasov-Poisson system in a half space. Preprint.

\bibitem {LionsPerthame}P. L. Lions, B. Perthame, Propagation of moments and
regularity of solutions for the 3-dimensional Vlasov-Poisson system. Invent.
Math. \textbf{105}, 415-430 (1991).

\bibitem {M}McOwen, R. C.: Partial Differential Equations, Methods and
Applications, Second Edition, Prentice Hall, 1996.

\bibitem {Pfaffelmoser}K. Pfaffelmoser, Global classical solutions of the
Vlasov-Poisson system in three dimensions for general initial data. J.
Differential Equations \textbf{95}, 281-303 (1992).

\bibitem {Sc1}J. Schaeffer, Global existence of smooth solutions to the
Vlasov-Poisson system in three dimensions, Comm. Partial Differential
Equations, 16(1991), pp. 1313--1335.

\bibitem {S}D. J. Struik, Lectures on Classical Differential Geometry, Dover
Publication, New York (1998).
\end{thebibliography}
\end{document}